\documentclass[a4paper]{amsart}

\usepackage{amssymb,latexsym,amsmath,graphicx}
\usepackage{eucal}
\makeatletter
\@addtoreset{equation}{section}\makeatother

\newtheorem{theo}{Theorem}[section]
\newtheorem{defi}[theo]{Definition}
\newtheorem{lem}[theo]{Lemma}
\newtheorem{prop}[theo]{Proposition}
\newtheorem{cor}[theo]{Corollary}

\newcommand{\mc}{\mathcal}
\newcommand{\rr}{\mathbb{R}}
\newcommand{\nn}{\mathbb{N}}
\newcommand{\cc}{\mathbb{C}}
\newcommand{\hh}{\mathbb{H}}
\newcommand{\zz}{\mathbb{Z}}

\newcommand{\la}{\lambda}
\newcommand{\eps}{\epsilon}

\newcommand{\pl}{\partial}
\newcommand{\x}{\times}

\newcommand{\til}{\widetilde}
\newcommand{\bbar}{\overline}

\newcommand{\ddens}{\Gamma_{0}^{\demi}}

\newcommand{\cjd}{\rangle}
\newcommand{\cjg}{\langle}

\newcommand{\demi}{\frac{1}{2}}
\newcommand{\ndemi}{\frac{n}{2}}
\newcommand{\tra}{\textrm{Tr}}
\newcommand{\beq}{\begin{equation}}
\newcommand{\eeq}{\end{equation}}

\newcommand{\grad}{\textrm{grad}}
\newcommand{\rang}{\textrm{rank}~}
\newcommand{\rangp}{\textrm{Rank}}
\newcommand{\trans}{{^t}\!}
\def\qed{\hfill$\square$}
\begin{document}
\title[Meromorphic properties of the resolvent]
{Meromorphic properties of the resolvent on
asymptotically hyperbolic manifolds}
\author[Colin Guillarmou]{Colin Guillarmou}
\address{Laboratoire de Math\'ematiques Jean Leray\\
         UMR 6629 CNRS/Universit\'e de Nantes \\
         2, rue de la Houssini\`ere \\
         BP 92208 \\
         44322 Nantes Cedex 03\\
         France}
     \email{cguillar@math.univ-nantes.fr}
%
\maketitle

\begin{abstract}
\noindent On an asymptotically hyperbolic manifold $(X^{n+1},g)$,
Mazzeo and Melrose have constructed the meromorphic extension
of the resolvent $R(\la):=(\Delta_g-\la(n-\la))^{-1}$ for the
Laplacian. However, there are special points on $\demi(n-\nn)$
that they did not deal with. We show that the points of
$\ndemi-\nn$ are at most some poles of finite multiplicity, and
that the same property holds for the points of $\frac{n+1}{2}-\nn$
if and  only if the metric is `even'. On the other hand, there
exist some metrics for which $R(\la)$ has an essential
singularity on $\frac{n+1}{2}-\nn$ and these cases are generic.
At last, to illustrate them, we give some examples with a sequence of
poles of $R(\la)$ approaching an essential singularity.
\end{abstract}
\vspace{0.5cm}

\section{Introduction and statements of the results}
The purpose of this work is to analyze, near the points
$\left(\frac{n-k}{2}\right)_{k\in\nn}$, the meromorphically
continued resolvent for the Laplacian
\[R(\la):=(\Delta_g-\la(n-\la))^{-1}\]
on some non-compact spaces $(X^{n+1},g)$ called
asymptotically hyperbolic manifolds.
This meromorphic extension in $\cc\setminus \demi(n-\nn)$
with finite rank poles, proved
by Mazzeo and Melrose \cite{MM}, is a beautiful application
of Melrose's pseudodifferential calculus on manifolds with
corners, which generalizes some well-known results on hyperbolic
spaces.
Meromorphic extensions of resolvents have been studied in many
frameworks and their finite rank poles, called resonances, serve
in a sense as discrete data similar in character to eigenvalues
of a compact manifold.
As far as we are concerned, the construction of \cite{MM} does not
treat the special points $\left(\frac{n-k}{2}\right)_{k\in\nn}$ and,
as Borthwick and Perry noticed in their article \cite{PB}, it seems
possible that these points are poles with infinite rank residues,
or even essential singularities of $R(\la)$.

However, if the manifold has constant negative sectional curvature
away from a compact, Guillop\'e and Zworski \cite{GZ2} did show
the meromorphic continuation of the resolvent to $\cc$ with finite rank poles.
The key to analyze the points of $\demi(n-\nn)$ is the special
structure of the metric near infinity, in the sense that its
Laplacian is locally the hyperbolic Laplacian, whose coefficients remain
smooth at $z=0$ in the new coordinates $(x_1,\dots,x_n,z=y^2)$ on $\hh^{n+1}$.
We could then follow the construction of \cite{MM} and search the
good conditions to set on the metric in order to use the same kind of
arguments: we would find that the natural assumption is to take
a metric with an even asymptotic expansion at infinity, in a sense we will
explain later.

In fact, our philosophy will be to use the properties of the scattering
operator, whose poles are essentially the resonances (cf. \cite{PB}).
The recent work of Graham and Zworski \cite{GRZ} gives indeed
a simple and explicit presentation of the scattering operator
$S(\la)$ on asymptotically hyperbolic manifolds which allows us
to study the nature of $S(\la)$ near $\demi(n+\nn)$. Thanks to
their calculus and the formula $S(n-\la)=S(\la)^{-1}$, we detail
the behavior of $S(\la)$ near
$\demi(n-\nn)$ and the relations between $R(\la)$ and $S(\la)$
provide a good analysis of the resolvent in a neighbourhood of the
points $\left(\frac{n-k}{2}\right)_{k\in\nn}$.\\

Firstly, let us recall some basic definitions and results to understand
the problem.
Let $\bar{X}=X\cup\pl\bar{X}$ a $n+1$-dimensional
smooth compact manifold with boundary and $x$ a defining function for
the boundary, that is a smooth function $x$ on $\bar{X}$ such that
\[x\geq 0,\quad \pl\bar{X}=\{m\in\bar{X}, x(m)=0\},\quad
dx|_{\pl\bar{X}}\not=0.\]
We say that a smooth metric $g$ on the interior $X$ of $\bar{X}$ is
\textsl{conformally compact} if $x^2g$ extends smoothly as a metric
to $\bar{X}$.
An \textsl{asymptotically hyperbolic manifold} is a
conformally compact manifold such that for all $y\in\pl\bar{X}$,
all sectional curvatures at $m\in X$ converge to $-1$ as $m\to y$.
Such a manifold is necessarily complete and the spectrum of its Laplacian
$\Delta_g$ acting on functions consists of absolutely continuous spectrum
$[\frac{n^2}{4},\infty)$ and a finite set of eigenvalues
$\sigma_{pp}(\Delta_g)\subset (0,\frac{n^2}{4})$.
The resolvent $(\Delta_g-z)^{-1}$ is a meromorphic family on
$\cc\setminus [\frac{n^2}{4},\infty)$ of bounded operators and
the new parameter $z=\la(n-\la)$ with $\Re(\la)>\ndemi$
induces a modified resolvent
\[R(\la):=(\Delta_g-\la(n-\la))^{-1}\]
which is meromorphic on $\{\Re(\la)>\ndemi\}$, his poles being
the points $\la_e$ such that $\la_e(n-\la_e)\in\sigma_{pp}(\Delta_g)$.
The problems studied by Mazzeo and Melrose in \cite{MM}
are the existence and the construction of the meromorphic extension
of $R(\la)$ to some open subsets of $\cc$ with values in weighted spaces.\\

Before we recall Mazzeo-Melrose theorem, let us introduce
a few notations to simplify the statements:
for $N\in\rr$ and $k=1,2$, let
\[\mc{O}_N:=\left\{\la\in\cc; \Re(\la)> \ndemi-N\right\},
\quad Z^k_{\pm}:=\ndemi\pm(\frac{k}{2}+\nn_0)\subset \cc.\] 
Let $(\mc{B}_i)_{i=0,1,2}$ be some Banach spaces, we denote
$\mc{L}(\mc{B}_1,\mc{B}_2)$ (or $\mc{L}(\mc{B}_1)$ if
$\mc{B}_1=\mc{B}_2$) the space of bounded linear operators from
$\mc{B}_1$ to $\mc{B}_2$. If $U$ is an open domain of $\cc$,
$\mc{H}ol(U,\mc{B}_1)$ (resp. $\mc{M}er(U,\mc{B}_1)$) is the set
of holomorphic (resp. meromorphic) functions on $U$ with values in
$\mc{B}_1$. $M(\la)$ is said meromorphic in $U$ with values in the
Banach space $\mc{B}_0$ if for each $\la_0\in U$, there exist a
neighbourhood $V_{\la_0}$ of $\la_0$, an integer $p>0$ and some
$(M_i)_{i=1,\dots,p}$ in $\mc{B}_0$ such that for all $\la\in
V_{\la_0}\setminus \{\la_0\}$ we have the finite Laurent expansion
\beq\label{serielaurent}
M(\la)=\sum_{i=1}^pM_i(\la-\la_0)^{-i}+H(\la), \quad
H(\la)\in\mc{H}ol(V_{\la_0},\mc{B}_0). 
\eeq
It is easy to see that $M(\la)$ is holomorphic in $U\setminus S$
where $S$ is a discrete set of $U$ whose elements are the poles
of $M(\la)$. $p$ is the order of the pole, $M_1=\textrm{Res}_{\la_0}M(\la)$ is the residue
of $M(\la)$ at $\la_0$ and if $\mc{B}_0=\mc{L}(\mc{B}_1,\mc{B}_2)$ is a
space of continuous linear maps, $m_{\la_0}(M(\la)):=\rang M_1$ is called the
multiplicity of $\la_0$ and $\rangp_{\la_0}M(\la):=\sum_{i=1}^p\rang M_i$
the total polar rank of $M(\la)$ at $\la_0$. If now the total polar rank is finite at
each pole of $M(\la)$ we say that $M(\la)$ is finite-meromorphic and
 $\mc{M}er_f(U,\mc{B}_0)$ denotes the space of finite-meromorphic
functions in $U$ with values in $\mc{B}_0$. At last,
if $\la_0\in U$ and $M(\la)$ is meromorphic in
$U\setminus\{\la_0\}$ but not in $U$, we will say that $\la_0$
is an essential singularity of $M(\la)$.

At last, note that all these definitions extend to locally convex vector
spaces (see for instance Bunke-Olbrich \cite{BO}).\\

Here is an interpretation of the result of Mazzeo and Melrose \cite[Th. 7.1]{MM} :
\begin{theo}\label{mazzeomelrose}
Let $(X,g)$ be an asymptotically hyperbolic manifold,
$\Delta_g$ its Laplacian acting on functions and $x$ a boundary defining function on
$\bar{X}$. The modified resolvent
\[R(\la):=(\Delta_g-\la(n-\la))^{-1} \in \mc{M}er_f \left(\mc{O}_0,
\mc{L}(L^2(X))\right)\]
with poles at points $\la\in \mc{O}_0$ such that
$\la(n-\la)\in\sigma_{pp}(P)$, extends to a finite-meromorphic family
\[R(\la)\in \mc{M}er_f\left(\mc{O}_N\setminus(Z^1_-\cup Z^2_-),
\mc{L}(x^NL^2(X),x^{-N}L^2(X))\right), \quad \forall N\geq 0\]
\end{theo}

The poles of this extension are called \textsl{resonances} and they
do not depend on $N$ (neither the multiplicity $m_{\la_0}((n-2\la)R(\la))$
of a resonance $\la_0$), they form a discrete set $\mc{R}\subset
\cc\setminus \{Z^1_-\cup Z^2_-\}$.
At each $\la\in(Z^1_-\cup Z^2_-)$, the behavior of $R(\la)$ is not clear,
it could be a pole of infinite multiplicity or an essential singularity.
Observe that the equation $(\Delta_g-\la(n-\la))R(\la)=1$ implies that
for a pole $\la_0$ of $R(\la)$, $\rangp_{\la_0}R(\la)<\infty$
is equivalent to $m_{\la_0}(R(\la))<\infty$. Let us now give a definition
which will be essential and will be explained in Section 2:
\begin{defi}\label{metricpaire}
Let $(X,g)$ be an asymptotically hyperbolic manifold and $k\in\nn\cup\{\infty\}$.
We say that $g$ is even modulo $O(x^{2k+1})$ if there exists $\eps>0$,
a boundary defining function $x$ and
some tensors $(h_{2i})_{i=0,\dots,k}$ on $\pl\bar{X}$ such that
\beq\label{metricpaire2}
\phi^*(x^2g)=dt^2+\sum_{i=0}^kh_{2i}t^{2i}+O(t^{2k+1})
\eeq
where $\phi$ is the diffeomorphism induced by the flow $\phi_t$
of the gradient $\grad_{x^2g}(x)$:
\[\phi :\left\{\begin{array}{rcl}
[0,\eps)\x\pl\bar{X}&\to &\phi([0,\eps)\x\pl\bar{X})\subset \bar{X}\\
(t,y)&\to&\phi_t(y)
\end{array}
\right.\]
\end{defi}

Using the relations between resolvent and scattering operator in a
way similar to \cite{CV,G1,PP} and the
calculus of the residues of $S(\la)$ by Graham-Zworski \cite{GRZ}
we find a necessary and sufficient condition on the metric to
have a finite-meromorphic extension of the resolvent to $\cc$.

\begin{prop}\label{amelioration}
Under the assumptions of Theorem \ref{mazzeomelrose}, the modified
resolvent extends to a finite-meromorphic family
\beq\label{prresolv}
R(\la)\in \mc{M}er_f\left(\mc{O}_N\setminus
Z^1_-, \mc{L}(x^NL^2(X),x^{-N}L^2(X))\right), \quad \forall N\geq 0
\eeq
and if $g$ is even modulo $O(x^{2k+1})$, this extension satisfies
\beq\label{prresolv2}
R(\la)\in\mc{M}er_f\left(\mc{O}_N,
\mc{L}(x^NL^2(X),x^{-N}L^2(X))\right), \quad \forall N\in[0,k+\demi)
\eeq
Conversely if (\ref{prresolv2}) holds true for $k\geq 2$ then
$g$ is even modulo $O(x^{2k-1})$.
\end{prop}

We then deduce the following
\begin{theo}\label{parite}
Let $(X,g)$ be an asymptotically hyperbolic manifold,
then $R(\la)$ admits a finite-meromorphic
extension to $\cc$ if and only if $g$ is even modulo $O(x^\infty)$.
\end{theo}

\textsl{Remark}: as a matter of fact, the usual examples
are some particular cases of
even metrics: the hyperbolic metrics perturbed on a compact
\cite{GZ1,GZ2,GZ3,PP}, the De
Sitter-Schwarzschild model \cite{SBZ},
the almost-product type metrics \cite{JSB}.
The asymptotically Einstein manifolds of dimension $n+1$
are only even modulo $O(x^n)$ in general \cite{GR,GRZ}.\\

Let us denote by $\mc{M}_{ah}(X)$ the space of asymptotically
hyperbolic metrics on $X$ with the topology inherited from
$x^{-2}C^\infty(\bar{X},T^*\bar{X}\otimes T^*\bar{X})$. If the
metric is not even, there is at least a point of $Z^1_-$ which is
either a pole of infinite multiplicity or an essential singularity
of $R(\la)$. The following results show that an essential
singularity appears generically.

\begin{theo}\label{genericitean}
Let $\bar{X}$ be a compact manifold with boundary of dimension
$n+1>2$. Then the set of metrics in $\mc{M}_{ah}(X)$ for which
$\frac{n-1}{2}$ is an essential singularity of $R(\la)$ contains
an open and dense set in $\mc{M}_{ah}(X)$.
\end{theo}

Note that the proof of Theorem \ref{genericitean} implies a little
more general result including the case $n=1$: the same is indeed
satisfied for the point $\frac{n-1}{2}-k$ if we consider the set of
even metrics modulo $O(x^{2k+1})$ instead of $\mc{M}_{ah}(X)$
and if $\frac{n-1}{2}-k\not=0$.

In dimension $n+1=2$, we will see that for $g\in\mc{M}_{ah}(X)$ analytic
near the boundary such that $\pl\bar{X}$ is a
connected geodesic of $(\bar{X},x^2g)$, the
resolvent is meromorphic if and only if $g$ is even.\\

Finally, to illustrate these results, we give some examples
with a sequence of resonances approaching
an essential singularity of $R(\la)$.

\begin{prop}\label{singessentielles}
For all $k\in\nn_0$ such that $2k\not=n-1$, there exists a
$n+1$-dimensional asymptotically hyperbolic manifold such that the
extension (\ref{prresolv}) has a sequence of poles which converges
to $\frac{n-1}{2}-k$.
\end{prop}

These cases are the first examples (as far as we know) of
essential singularities coming from the meromorphic extension
of the resolvent. If the resonances are interpreted as
eigenvalues of an operator (see \cite{A} when the extension is
finite-meromorphic), we could
think that these essential singularities are some
isolated points in the essential spectrum of this operator.

\begin{figure}[ht!]
\begin{center}
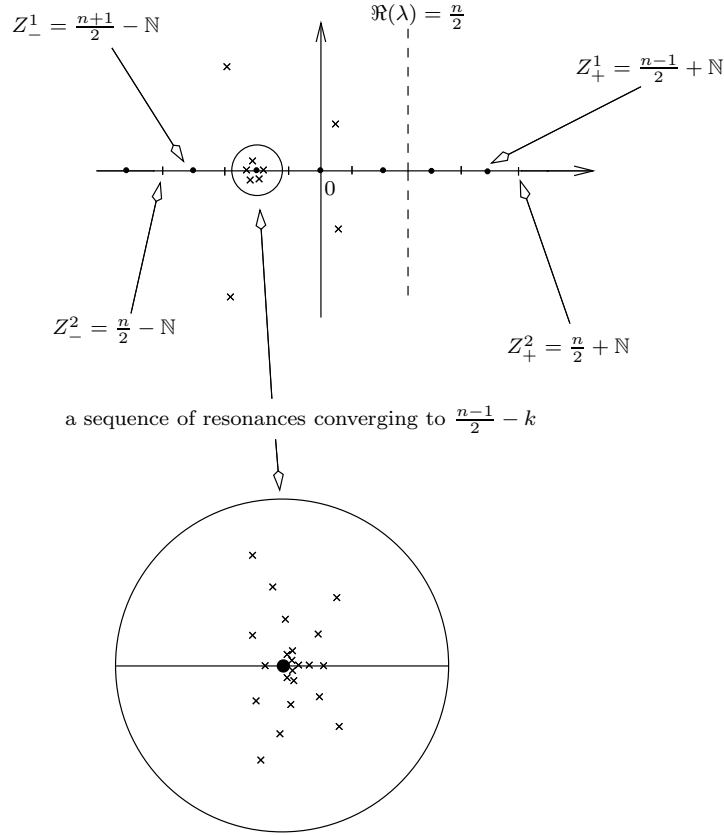
\caption{The resonances of $\Delta_g$ in a case where $g$ is not even}
\end{center}
\end{figure}

The paper is organized as follows:
Section 2 recalls some basic geometric facts on asymptotically hyperbolic
manifolds and explains Definition \ref{metricpaire};
Section 3 shows how to use the scattering operator instead of the resolvent;
we give in Section 4 the proofs of the main results and Section 5 contains
the examples.\\

\textbf{Acknowledgements.} I would like to thank L. Guillop\'e, G. Carron, D. Borthwick,
G. Vodev and N. Yeganefar for their help and comments. I also
thank R. Graham for a discussion about asymptotically hyperbolic
Einstein metrics.

\section{Geometry of $(X,g)$, even metric}

Let $(X,g)$ an asymptotically hyperbolic manifold,
$x_0$ a boundary defining function and
\[H:=x_0^2g\]
the induced metric by $g$ and $x_0$ on $\bar{X}$.
We can easily check that neither the conformal
class $[H|_{T\pl\bar{X}}]$
of the metric $H|_{T\pl\bar{X}}$ on
$\pl\bar{X}$, nor the value at the boundary
$(|dx_0|_{H})|_{\pl\bar{X}}$ depends
on the choice of the function $x_0$.
Moreover, Mazzeo and Melrose \cite{MM} remark that the
sectional curvatures of $g$ at $m\in X$ approach
$-|dx_0|^2_H(y)$ when $m\to y\in\pl\bar{X}$, so we can summarize the
property `asymptotically hyperbolic' with the identity
$|dx_0|_{H}=1$ on $\pl\bar{X}$ (which does not depend on the choice of $x_0$).

If $(X,g)$ is asymptotically hyperbolic, it is shown by
Graham \cite{GR} that there exists, for each metric
$h_0\in[H|_{T\pl\bar{X}}]$, a unique
boundary defining function $x$ such that $|dx|_{x^2g}=1$ in an
open neighbourhood $V_x \subset \bar{X}$ of $\pl\bar{X}$
and $(x^2g)|_{T\pl\bar{X}}=h_0$.
So there is a collar $U_x:=[0,\eps_x)\x\pl\bar{X}$
linked to $x$ by the diffeomorphism
\beq\label{collier}
\phi :\left\{\begin{array}{rcl}
U_x&\to &\phi(U_x)\subset V_x\\
(t,y)&\to&\phi_t(y)
\end{array}
\right.
\eeq
where $\phi_t$ is the flow of the gradient $\grad_{x^2g}(x)$.
Note that $t=x$ as functions on $\phi(U_x)$.
In the open collar $(0,\eps_x)\x\pl\bar{X}$, the metric $g$
can be expressed by
\beq\label{formenormale}
\phi^* g=\frac{dt^2+h(t,y,dy)}{t^2},\quad h(0,y,dy)=h_0(y,dy)
,\quad h\in C^\infty(U_x,S^2(T^*U_x))\eeq
$S^2(T^*U_x)\subset T^*U_x\otimes T^*U_x$ being
the bundle over $U_x$ of symmetric 2-tensors.
This form is called a \textsl{model form} and we shall
write $g$ and $x$ instead of $\phi^*g$ and $t$ in
(\ref{formenormale}) for simplicity.
Let us now define the set of boundary defining functions that
induce a model form of the metric:
\[Z(\pl\bar{X}):=\{x\in C^{\infty}(\bar{X});x\geq 0,
\pl\bar{X}=x^{-1}(0), \exists \eps_x>0,\forall m\in x^{-1}([0,\eps_x)),
|dx|_{x^2g}(m)=1\}\]
for which Graham \cite{GR} has constructed a bijection
\[ [H|_{T\pl\bar{X}}]\longleftrightarrow Z(\pl\bar{X}).\]
The symmetric tensor $h(t,y,dy)$ in (\ref{formenormale})
defines a family of metrics $h(t)$ on the hypersurfaces
$\{x=t\}$, and it depends on the choice of the function
$x\in Z(\pl\bar{X})$.
We can easily check that, for a fixed $k\in\nn$,
the vanishing condition modulo $O(x^k)$ at $\pl\bar{X}$
\[\quad h(x)-h(0)=O(x^k)\]
is invariant with respect to the boundary defining function $x\in Z(\pl\bar{X})$.
But to treat our problem it is more natural to choose the weaker condition
introduced in Definition \ref{metricpaire}: there exists $x\in Z(\pl\bar{X})$
and $k\in\nn$ such that the Taylor expansion
of $x^2g$ at $x=0$ consists only of even powers of $x$ up through
the $x^{2k+1}$ term
\beq\label{metriquepaire}
x^2g=dx^2+l(x^2,y,dy)+O(x^{2k+1}), \quad l\in C^{\infty}(U_x,S^2(T^*U_x))
\eeq
in the collar $U_x$ linked to $x$ by (\ref{collier}).
As a matter of fact, this property is not associated to a particular
defining function, as we could think, but to the set
$Z(\pl\bar{X})$ or equivalently to the conformal class
$[H|_{T\pl\bar{X}}]$: indeed, if the property is satisfied for one
function of $Z(\pl\bar{X})$, it is satisfied for all functions of
$Z(\pl\bar{X})$. The metric is then
said to be \textsl{even modulo $O(x^{2k+1})$}.

\begin{lem}\label{invarianceparite}
Let $(X,g)$ be an asymptotically hyperbolic manifold. Suppose that
there exists a function $x\in Z(\pl\bar{X})$ and $k\in\nn$ such that
the metric $x^2g$ can be expressed by
\beq\label{metriquepaire2}
x^2g=dx^2+l(x^2,y,dy)+O(x^{2k+1}), \quad
l\in C^{\infty}(U_x,S^2(T^*U_x))
\eeq
in the collar $U_x=[0,\eps_x)\x\pl\bar{X}$ linked to $x$ by (\ref{collier}).
Then, for all function $t\in Z(\pl\bar{X})$, the metric $t^2g$
can be expressed by
\[t^2g=dt^2+p(t^2,z,dz)+O(t^{2k+1}), \quad
p\in C^{\infty}(U_t,S^2(T^*U_t))\]
in the collar $U_t=[0,\eps_t)\x\pl\bar{X}$ linked to $t$ by (\ref{collier}).
\end{lem}
\textsl{Proof}: firstly, we recall Graham calculus in \cite[Lem. 2.1 and 2.2]{GR}.
Let $x\in Z(\pl\bar{X})$ such that (\ref{metriquepaire2}) holds and
$t\in Z(\pl\bar{X})$. Let $h$ defined as in (\ref{formenormale}) and set
\[t=e^{\omega}x, \quad \omega\in C^\infty(U_x)\]
According to \cite{GR}, $\omega$ is a solution of the non-linear equation
\beq\label{eqnonlin}
2\pl_x\omega+x\Big((\pl_x\omega)^2+\sum_{ij}h^{ij}(x)
\pl_{y_i}\omega\pl_{y_j}\omega\Big)=0.
\eeq
We easily obtain $\pl_x\omega|_{x=0}=0$ and by differentiating (\ref{eqnonlin})
an even number of times with respect to $x$, it can be shown by induction
that $\pl_x^{2j+1}\omega|_{x=0}=0$ for $j\leq k$,
(cf. \cite{GR} for details).
Recall now that the collar linked to $t$ is constructed by the diffeomorphism
\[\phi' :\left\{\begin{array}{rcl}
U_t:=[0,\eps_t)\x\pl\bar{X}&\to &\phi'(U_t)\\
(t,z)&\to&\phi'_t(z)
\end{array}
\right.\]
where $\phi'_t(z)$ is the flow of $\grad_{t^2g}t$. Let us set
\[x(t,z):=x(\phi'_t(z)),\quad y(t,z):=y(\phi'_t(z))\]
and we will show by induction that we have for all $m\leq 2k+2$
\beq\label{recurparite}
\begin{array}{ll}
\pl_t^{2j}x(t,z)|_{t=0}=0 & \forall j, 0\leq 2j\leq m\\
\pl_t^{2j+1}y(t,z)|_{t=0}=0 & \forall j, 0\leq 2j+1\leq m.
\end{array}
\eeq
To begin, note that
\begin{eqnarray*}
\grad_{t^2g}t&=&t^{-2}\grad_{g}t=e^{-\omega}\grad_{x^2g}x+
e^{-\omega}x\grad_{x^2g}\omega\\
&=& (e^{-\omega}+te^{-2\omega}\pl_x\omega)\pl_x+e^{-2\omega}t\sum_{i,j}h^{ij}(x)
\pl_{y_i}\omega\pl_{y_j}
\end{eqnarray*}
and $x(t,z)$, $y(t,z)$ are defined by the following flow equations
\begin{eqnarray}
\label{eqnflot1}
\pl_t x(t,z)&=&e^{-\omega(x(t,z),y(t,z))}+
te^{-2\omega (x(t,z),y(t,z))}\pl_x\omega(x(t,z),y(t,z)),\\
\label{eqnflot2}
\pl_{t}y_j(t,z)&=&e^{-2\omega(x(t,z),y(t,z))}t
\sum_{i=1}^nh^{ij}\left(x(t,z),y(t,z)\right)
\pl_{y_i}\omega(x(t,z),y(t,z)).
\end{eqnarray}
To show (\ref{recurparite}), we first check it for $m=1$
\[x(0,z)=0,\quad  \pl_t y_j(0,z)=0,\quad j=1,\dots,n. \]
Suppose now that (\ref{recurparite}) is satisfied for the integer $m$
(with $m\leq 2k+1$).

If $m+1$ is even and  $m+1\leq 2k+2$, we obtain from (\ref{eqnflot1})
\[\pl^{m+1}_tx(0,z)=\pl^m_te^{-\omega}|_{t=0}+
m\pl^{m-1}_t(e^{-2\omega}\pl_x\omega)|_{t=0}.\]
Remark that if $f(x,y)$ is an arbitrary even function (resp. odd)
in $x$ modulo $O(x^{2l+1})$ (resp. modulo $O(x^{2l})$),
then the composition of $f(x,y)$ with
\[(t,z)\to(x(t,z),y(t,z))\]
is even (resp. odd) in $t$ modulo $O(t^{\min(2l+1,m+2)})$
(resp. modulo $O(t^{\min(2l,m+1)})$).
Therefore, we deduce that
\[\pl^m_te^{-\omega}|_{t=0}=0.\]
because $m$ is odd and $e^{-\omega}$ is even modulo $O(x^{2k+3})$.
Moreover, $m-1$ is even and the derivatives
$\pl^{m-1}_t(e^{-2\omega}\pl_x\omega)|_{t=0}$
split into a sum of products of derivatives of $e^{-2\omega}$ and
$\pl_x\omega$. In each product, if we differentiate
an odd number of times one of the terms
at $t=0$, the number of derivatives in the other term must be odd
too and the previous argument shows that the product vanishes because
$e^{-2\omega}$ is even in $x$ modulo $O(x^{2k+3})$.
If the number of derivatives for one term of the product is even,
the number for the other term is even too and the product vanishes
because $\pl_x\omega$ is odd in $x$ modulo $O(x^{2k+2})$.
We then deduce that
\[\pl^{m+1}_tx(0,z)=0\]

On the other hand, if $m+1$ is odd and $m+1\leq 2k+1$, we use
the same trick for equation
(\ref{eqnflot2}) and we just have to differentiate an odd ($=m-1$) number
of times a product of even functions in $t$
modulo $O(t^{\min(2k+1,m+1)})$, which proves that
\[\pl^{m+1}_ty_i(0,z)=0\]
and we conclude by induction that (\ref{recurparite}) is true for all $m\leq 2k+2$.

We finally have to show that for all $\xi\in T_z{\pl\bar{X}}$
\[t^2g(\xi,\xi)=e^{2\omega}\left((\pl_zx.dz(\xi))^2+
h(x,y,\pl_zy.dz(\xi))\right)\]
is an even function in $t$ modulo $O(t^{2k+1})$, which is a simple consequence
of the odd-even properties of $\omega$, $x$, $y$ and $h$.
\qed\\

Let us denote by $\bar{X}^2:=(\bar{X}\bigsqcup\bar{X})/\pl\bar{X}$
the double of $\bar{X}$, which is firstly a topological space.
Choose $x$ a boundary defining function of $\pl\bar{X}$.
From the diffeomorphism (\ref{collier}),
we can construct a $C^\infty$ atlas
on $\bar{X}^2$, using the fact that $\pl\bar{X}\subset \bar{X}^2$ is
contained in a an open set $V^2_x:=(V_1\bigsqcup V_2)/\pl\bar{X}$
(with $V_1=V_2=\phi(U_x)$) diffeomorphic to $(-\eps_x,\eps_x)\x\pl\bar{X}$ via
\begin{eqnarray*}
 (-\eps_x,\eps_x)\x \pl\bar{X}&\simeq & V^2_x\\
 (t,y)&\to &\left\{
\begin{array}{ll}
 [\phi_{-t}(y)], & \phi_{-t}(y)\in V_1\textrm{ if } t\leq0\\
 \lbrack\phi_t(y)\rbrack, & \phi_t(y)\in V_2\textrm{ if } t\geq0
\end{array}\right..
\end{eqnarray*}
The remaining charts come easily from the charts of the interior
$X$ of $\bar{X}$.
Remark that this $C^\infty$ structure on
$\bar{X}^2$ depends on the choice of $x$.
If now we denote this structure by $\bar{X}^2_x$,
there is a global diffeomorphism
\[\bar{X}^2_x\simeq \bar{X}^2_{x'}\]
for two different boundary defining functions $x$ and $x'$,
but it is not the case that any one of these structures is natural with
respect to $C^\infty(\bar{X})$.
By the even functions of $x$ modulo $O(x^{2k+1})$
on $\bar{X}$, we shall mean the smooth functions on $\bar{X}$ which
admit a $C^{2k}$ continuation to $\bar{X}^2_x$ and are
invariant with respect to the natural involution
exchanging the factors on $\bar{X}^2$.
The result is a class of functions which depends on the choice of $x$:
a function whose Taylor expansion in $x$ at $x=0$ is even
modulo $O(x^{2k+1})$ does not necessarily have an even
Taylor expansion in $x'$ modulo $O(x'^{2k+1})$ at $x'=0$, if $x,x'$ are
two different boundary defining functions of $\bar{X}$.
In the proof of Lemma \ref{invarianceparite}, it is shown
that if the metric can be expressed by (\ref{metriquepaire}) for one
boundary defining function $x$, then the coordinate changes
$(x,y)\to(x',y')$ on $[0,\eps)\x\pl\bar{X}$
which leave the metric under a model form have local expansions of the form
\[x'=x\sum_{j=0}^{k+1}a_j(y)x^{2j}+O(x^{2k+4}),\quad
y'=\sum_{j=0}^{k+1}b_j(y)x^{2j}+O(x^{2k+3})\]
with some smooth functions $a_j$ and $ b_j$,
thus they induce some $C^{2k+2}$ compatible charts on $\bar{X}^2_x$.
As a conclusion, if $x\in Z(\pl\bar{X})$ the structure $\bar{X}^2_x$
does not depend on the choice of $x$ as a $C^{2k+2}$
structure on $\bar{X}^2$ and we obtain a natural choice (with respect to $g$) of
$C^{2k+2}$ structure on $\bar{X}^2$ induced by the functions
in $Z(\pl\bar{X})$. Moreover, it admits a $C^{2k}$ conformal class
of metrics which are invariant with respect to
the natural involution exchanging the factors on $\bar{X}^2$: this is
obtained by extending by symmetry the metrics $x^2g$
for each $x\in Z(\pl\bar{X})$.

\section{From the resolvent to the scattering operator}

\subsection{Stretched products}

To begin, let us introduce a few notations and recall
some basic things on stretched products (the reader can
refer to Mazzeo-Melrose \cite{MM}, Mazzeo \cite{MA1}
or Melrose \cite{M1} for details).
Let $\bar{X}$ a smooth compact manifold with boundary and $x$
a boundary defining function.
The manifold $\bar{X}\x\bar{X}$ is a smooth manifold with corners,
whose boundary hypersurfaces are diffeomorphic to
$\pl\bar{X}\x\bar{X}$ and $\bar{X}\x\pl\bar{X}$, and defined by
the functions $\pi_L^*x$, $\pi_R^*x$ ($\pi_L$ and $\pi_R$ being the
left and right projections from $\bar{X}\x\bar{X}$ onto $\bar{X}$).
For notational simplicity, we now write $x$,$x'$ instead of $\pi_L^*x$, $\pi_R^*x$
and let
\[\delta_{\pl\bar{X}}:=\{(m,m)\in\pl\bar{X}\x\pl\bar{X}; m\in \pl\bar{X}\}.\]
Remark the elementary embeddings
\beq\label{plongements}
\delta_{\pl\bar{X}}\to \pl\bar{X}\x\pl\bar{X}\to\pl\bar{X}\x\bar{X}\to\bar{X}\x\bar{X}
\eeq
which will often be regarded as inclusions.
The blow up of $\bar{X}\x\bar{X}$ along the `diagonal' $\delta_{\pl\bar{X}}$
of $\pl\bar{X}\x\pl\bar{X}$ will be denoted by $\bar{X}\x_0\bar{X}$
and the blow-down map
\[\beta:\bar{X}\x_0\bar{X}\to \bar{X}\x\bar{X}.\]
This manifold with corners has three boundary hypersurfaces
$\mc{T},\mc{B},\mc{F}$ defined by some functions $\rho,\rho',R$
such that $\beta^*(x)=R\rho$, $\beta^*(x')=R\rho'$.
Globally, $\delta_{\pl\bar{X}}$ is replaced by a larger manifold, namely
by its doubly inward-pointing spherical normal bundle of $\delta_{\pl\bar{X}}$,
whose each fiber is a quarter of sphere. From local coordinates
$(x,y,x',y')$ on $\bar{X}\x\bar{X}$, this amounts to introducing
polar coordinates  $(R,\rho,\rho',\omega,y)$
around $\delta_{\pl\bar{X}}$:
\[R:=(x^2+x'^2+|y-y'|^2)^\demi, \quad (\rho,\rho',\omega):=
\left(\frac{x}{R},\frac{x'}{R},\frac{y-y'}{R}\right)\]
with $R,\rho,\rho'\in [0,\infty)$. These polar coordinates are useful to
describe the singularities of the Schwartz kernel of $R(\la)$.

\begin{figure}[ht!]
\begin{center}
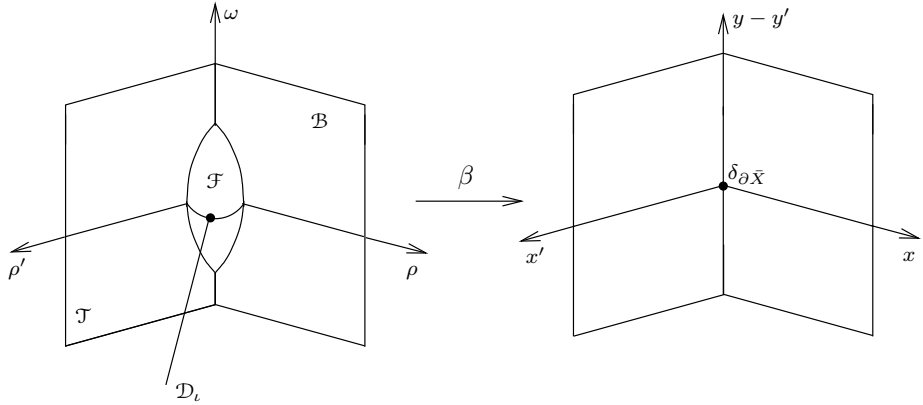
\caption{The blow-down map}\label{blow}
\end{center}
\end{figure}

Similarly, we denote by $\pl\bar{X}\x_0\bar{X}$ the
blow-up of $\pl \bar{X}\x\bar{X}$ along $\delta_{\pl\bar{X}}$. It can be
naturally embedded in $\bar{X}\x_0\bar{X}$ with respect to
(\ref{plongements}):
\[\pl\bar{X}\x_0\bar{X}\simeq (\bar{X}\x_0\bar{X})\cap \mc{T}\]
using $\beta(\mc{T})=\{x=0\}\simeq \pl\bar{X}\x\bar{X}$.
With these identifications, $\til{\beta}:=\beta|_{\mc{T}}$ is
the blow-down map
\[\til{\beta}:\pl\bar{X}\x_0\bar{X}\to\pl\bar{X}\x\bar{X}.\]
This manifold with corners has two boundary hypersurfaces
$\til{\mc{B}},\til{\mc{F}}$ defined by the functions
$\til{\rho}':=\rho'|_{\mc{T}}$ and $\til{R}:=R|_{\mc{T}}$.\\

Finally, the blow-up $\pl\bar{X}\x_0\pl\bar{X}$ of $\pl\bar{X}\x\pl\bar{X}$
along $\delta_{\pl\bar{X}}$ can be naturally embedded
in $\pl\bar{X}\x_0\bar{X}$ with respect to
(\ref{plongements})
\[\pl\bar{X}\x_0\pl\bar{X}\simeq (\pl\bar{X}\x_0\bar{X})\cap \til{\mc{B}}\]
and the blow-down map is
$\beta_{\pl}:=\til{\beta}|_{\til{\mc{B}}}$. The function
$r:=\til{R}|_{\til{\mc{B}}}$ defines the boundary of
$\pl\bar{X}\x_0\pl\bar{X}$, which is the lift of $\delta_{\pl\bar{X}}$
under $\beta_\pl$. Let
$(y_0,y_0)\in\delta_{\pl\bar{X}}$, $V_{y_0}$ an open neighbourhood of
this point and $(y,y')$ a coordinate patch in $V_{y_0}$, then
\[\beta_\pl^*(|y-y'|)=rf,\quad f\in C^\infty(\beta_\pl^{-1}(V_{y_0})),\quad
f>0, \quad \pl_rf|_{r=0}\not=0.\]

\subsection{Half densities}
Let $\Gamma_0^\demi(\bar{X})$ the line bundle of singular
half-densities on $\bar{X}$, trivialized by $\nu:=|dvol_g|^{\demi}$, and
$\Gamma^\demi(\pl\bar{X})$ the bundle of half densities on $\pl\bar{X}$,
trivialized by $\nu_0:=|dvol_{h_0}|^\demi$ (where $h_0=x^2g|_{T\pl\bar{X}}$).
From these bundles, one can construct the bundles
$\Gamma_0^\demi(\bar{X}\x\bar{X})$, $\Gamma_0^\demi(\pl\bar{X}\x\bar{X})$
and $\Gamma^\demi(\pl\bar{X}\x\pl\bar{X})$ by tensor products,
whose sections are respectively
$\nu\otimes\nu$, $\nu_0\otimes\nu$ and $\nu_0\otimes\nu_0$.
Finally, let
$\Gamma_0^\demi(\bar{X}\x_0\bar{X})$ ,
$\Gamma_0^\demi(\pl\bar{X}\x_0\bar{X})$ and
$\Gamma^\demi(\pl\bar{X}\x_0\pl\bar{X})$ the bundles obtained by
lifting under $\beta$, $\til{\beta}$ and $\beta_{\pl}$ the three previous
bundles.
If $M$ denotes $\bar{X}$, $\bar{X}\x\bar{X}$ or $\pl\bar{X}\x\bar{X}$,
we write $\dot{C}^\infty(M,\Gamma_0^\demi(M))$ the space of smooths sections
of $\Gamma_0^\demi(M)$ that vanish to all order at all the boundary
hypersurfaces of $M$, and
$C^{-\infty}(M,\Gamma_0^\demi(M))$ is its topological dual, whose elements
are the extendible distributional half densities.
The Hilbert space $L^2(\bar{X},\Gamma_0^\demi(\bar{X}))$ constructed by
completing $\dot{C}^\infty(\bar{X},\Gamma_0^\demi(\bar{X}))$
with respect to the norm
\[||f||_2=\left(\int_{\bar{X}}f.\bar{f}\right)^\demi\]
is isomorphic to $L^2(X,dvol_g)$ and will be denoted by $L^2(X)$. Similarly it will
be more practical to write $H^R(\pl\bar{X})$ for the Sobolev space of order $R$
on $\pl\bar{X}$ with values in half densities. At last,
$\Gamma_0^\demi(M)$ and $\Gamma^\demi(M)$ will be replaced by
$\Gamma_0^\demi$ and $\Gamma^\demi$,
where $M$ is one of the previously introduced manifolds.\\

Set $(.,.)$ the symmetric non-degenerate products on
$L^2(X)$ and $L^2(\pl\bar{X})$
\[ (u,v):=\int_{X}uv, \quad (u,v):=\int_{\pl\bar{X}}uv.\]
For $\alpha\in\rr$, we can check by using the first product that the
dual space of $x^\alpha L^2(X)$ is isomorphic to
$x^{-\alpha}L^2(X)$. We shall also use the following tensorial
notation for $E=x^{\alpha} L^2(X)$ (resp. $E=L^2(\pl\bar{X})$),
$\psi,\phi\in E'$
\[\phi\otimes\psi : \left\{
\begin{array}{rcl}
E & \to & E'\\
f&\to& \phi( \psi,f)
\end{array}\right..\]

\subsection{Resolvent}
The meromorphically continued resolvent $R(\la)$
on $\cc\setminus (Z^1_-\cup Z^2_-)$
is a continuous operator from $\dot{C}^\infty(\bar{X},\Gamma_0^\demi)$
to $C^{-\infty}(\bar{X},\Gamma_0^\demi)$, its associated Schwartz kernel being
\[r(\la)\in C^{-\infty}(\bar{X}\x \bar{X},
\Gamma_0^\demi)\]
whose properties, studied in \cite{MM},
are recalled for instance in \cite[Th. 2.1]{PB}, namely
\[r(\la)=r_0(\la)+r_1(\la)+r_2(\la)\]
\[\beta^*(r_0(\la))\in I^{-2}(\bar{X}\x_0\bar{X},\Gamma_0^\demi)\]
\beq\label{k1}
\beta^*(r_1(\la))\in \rho^\la\rho'^\la C^\infty
(\bar{X}\x_0\bar{X},\Gamma_0^\demi)
\eeq
\beq\label{k2}
r_2(\la)\in x^{\la}x'^{\la}C^\infty(\bar{X}\x\bar{X},
\Gamma_0^\demi)
\eeq
where $I^{-2}(\bar{X}\x_0\bar{X},\Gamma_0^\demi)$ denotes
the set of conormal distributions of degree $-2$ on $\bar{X}\x_0\bar{X}$
associated to the closure of the lifted interior diagonal (see
Fig. \ref{blow})
\[\mc{D}_\iota:=\bbar{\beta^{-1}(\{(m,m)\in\bar{X}\x\bar{X}; m\in X\})}\]
and vanishing to infinite order at $\mc{B}\cup\mc{T}$ (notice that the
lifted interior diagonal only intersects the topological
boundary of $\bar{X}\x_0\bar{X}$ at $\mc{F}$ and it does transversally).
Moreover, $(\rho\rho')^{-\la}\beta^{*}(r_1(\la))$ and
$(xx')^{-\la}r_2(\la)$ are respectively holomorphic and meromorphic
in $\la\in\cc\setminus (Z^1_-\cup Z^2_-)$ and $r_0(\la)$ is the kernel
of a holomorphic family of operators
\[R_0(\la)\in \mc{H}ol(\cc, \mc{L}
(x^{\alpha}L^2(X),x^{-\alpha}L^2(X))),\quad
\forall \alpha\geq 0\]
Note also that Patterson-Perry arguments \cite[Lem. 4.9]{PP} prove that
$R(\la)$ does not have poles on the line $\{\Re(\la)=\ndemi\}$, except maybe
$\la=\ndemi$, it is a consequence of the absence
of embedded eigenvalues (see Mazzeo \cite{MA}). The only poles of $R(\la)$ in the half plane $\{\Re(\la)>\ndemi\}$
is the finite set of $\la_e$ such that
$\la_e(n-\la_e)\in \sigma_{pp}(\Delta_g)$, they are first order poles
and their residue is
\beq\label{residuvp}
\textrm{Res}_{\la_e}R(\la)=(2\la_e-n)^{-1}\sum_{k=1}^p \phi_k\otimes\phi_k,
\quad \phi_k\in x^{\la_e}C^\infty(\bar{X},\Gamma_0^\demi)
\eeq
where $(\phi_k)_{k=1,\dots,p}$ are the normalized
eigenfunctions of $\Delta_g$ for the eigenvalue $\la_e(n-\la_e)$.\\

\subsection{Poisson operator}
For $R(\la)\geq \ndemi$, $\la\notin \demi(n+\nn)\cup\mc{R}$
and $x$ a fixed boundary defining function of $\bar{X}$,
the Poisson operator is the unique continuous operator
\[\mc{P}(\la):C^\infty(\pl\bar{X},\Gamma^\demi)
\to x^{\la}C^\infty(\bar{X},\Gamma_0^\demi)
+x^{n-\la}C^\infty(\bar{X},\Gamma_0^\demi)\]
such that
\beq\label{problempoisson}
\left\{\begin{array}{l}
(\Delta_g-\la(n-\la))\mc{P}(\la)=0\\
\mc{P}(\la)f=x^{n-\la}F_1(\la)+x^{\la}F_2(\la)\\
F_1(\la),F_2(\la)\in C^\infty(\bar{X},\Gamma_0^\demi)\\
(x^{\ndemi}F_1(\la))|_{\pl\bar{X}}=f
\end{array}\right.
\eeq
Graham and Zworski \cite{GRZ} gave a
simple construction of $\mc{P}(\la)$ and
Joshi-S\'a Barreto \cite{JSB} proved that
its Schwartz kernel is a `weighted restriction'
of the resolvent kernel $r(\la)$ on the boundary $\{x'=0\}$, which implies
that $\mc{P}(\la)$ can be meromorphically continued to
$\cc\setminus (Z^1_-\cup Z^2_-)$.

For what follows, we now define the operator $E(\la)$ whose Schwartz kernel
$e(\la)$ is the weighted restriction of $r(\la)$ on the boundary $\{x=0\}$:
\beq\label{eisensteinfct}
e(\la):=\til{\beta}_*\left (\beta^*(x^{-\la+\ndemi}r(\la))|_{\mc{T}}\right )\in
C^{-\infty}(\pl\bar{X}\x \bar{X},\Gamma_0^\demi).
\eeq
According to (\ref{k1}) and (\ref{k2}),
the distribution $e(\la)$ has conormal
singularities on the boundaries described by
\[\til{\beta}^*(e(\la))\in
\til{\rho}'^{\la}\til{R}^{-\la+\ndemi}C^{\infty} (\pl\bar{X}\x_0
\bar{X}, \Gamma_0^\demi)+\til{\beta}^*(x'^\la
C^\infty(\pl\bar{X}\x\bar{X}, \Gamma_0^\demi)).\]
We then easily deduce that $E(\la)$ is continuous from
$\dot{C}^\infty(\bar{X},\Gamma_0^\demi)$ to
$C^\infty(\pl\bar{X},\Gamma^\demi)$, and its tranpose is well
defined from $C^{-\infty}(\pl\bar{X},\Gamma^\demi)$ to
$C^{-\infty}(\bar{X},\Gamma_0^\demi)$.
As a matter of fact, it follows from Joshi-S\'a Barreto work
\cite{JSB} that
\[\mc{P}(\la)=\trans E(\la)\]
on $C^{\infty}(\pl\bar{X},\Gamma^\demi)$.
Following the holomorphic-meromorphic properties of (\ref{k1}) and
(\ref{k2}), we obtain
\[\til{\beta}^*(x'^Ne(\la))\in \mc{H}ol (U_N,L^2(\pl\bar{X}\x_0\bar{X},\Gamma_0^\demi))\]
\[U_N:=\{\la\in\cc; \ndemi-N<\Re(\la)< N, \la\notin (\mc{R}\cup Z^1_-\cup Z^2_-)\}\]
and using the isomorphism $\til{\beta}^*$ between
$L^2(\pl\bar{X}\x\bar{X},\Gamma_0^\demi)$
and $L^2(\pl\bar{X}\x_0\bar{X},\Gamma_0^\demi)$,
we see that $E(\la)$ is Hilbert-Schmidt from
$x^NL^2(X)$ to $L^2(\pl\bar{X})$, and
\beq \label{propelambda}
E(\la)\in \mc{H}ol\left(U_N,\mc{L}(x^NL^2(X),L^2(\pl\bar{X}))\right).
\eeq
It also implies that its tranpose $\trans E(\la)$ is
holomorphic in $U_N$ with values in $\mc{L}(L^2(\pl\bar{X}),x^{-N}L^2(X))$.
Finally, the only  poles of $E(\la)$ and $\trans E(\la)$
in $\{\Re(\la)>\ndemi\}$ are the complex numbers
$\la_e$ such that $\la_e(n-\la_e)\in\sigma_{pp}(\Delta_g)$, (\ref{residuvp}) proving
that they are of order $1$ with finite multiplicity.\\

\subsection{Scattering operator}
For $\Re(\la)\geq \ndemi$, $\la\notin \demi(n+\nn)\cup\mc{R}$,
the scattering operator $S(\la)$ is defined by
\[S(\la):\left\{\begin{array}{ccc}
C^{\infty}(\pl\bar{X},\Gamma^\demi)&\to &C^{\infty}(\pl\bar{X},\Gamma^\demi)\\
f&\to &(x^\ndemi F_2(\la))|_{\pl\bar{X}}
\end{array}\right.\]
using the notations (\ref{problempoisson}). Using
the meromorphic continuation of $R(\la)$ to $\cc\setminus (Z^1_-\cup Z_-^2)$,
Joshi-S\'a Barreto \cite{JSB} proved that
$S(\la)$ can be meromorphically continued (weakly) to the same
set with Schwartz kernel $s(\la)$
\beq\label{noyaudiffusion}
s(\la):=(2\la-n)\left (\beta_{\pl}\right )_*\left (\beta^*\left (
x^{-\la+\ndemi}x'^{-\la+\ndemi}r(\la)\right )|_{\mc{T}\cap\mc{B}}\right )
\eeq
that is, in view of (\ref{k2}) and (\ref{noyaudiffusion}),
\beq\label{noyaudesla}
s(\la)=\left (\beta_{\pl}\right )_*\left (r^{-2\la}k_1(\la)\right )+
k_2(\la)
\eeq
\[k_1(\la):=(2\la-n)((\rho\rho')^{-\la+\ndemi}R^nr_1(\la))|_
{\mc{T}\cap\mc{B}}\in C^\infty(\pl\bar{X}\x_0\pl\bar{X},
\Gamma^\demi)\]
\[k_2(\la):=(2\la-n)((xx')^{-\la+\ndemi}r_2(\la))|_{x=x'=0}\in
C^\infty(\pl\bar{X}\x\pl\bar{X},\Gamma^\demi)\]
with $k_1(\la)$ and $k_2(\la)$ respectively holomorphic and meromorphic in
$\la\in \cc\setminus (Z^1_-\cup Z^2_-)$.
To understand the distribution $r^{-2\la}k_1(\la)$ on
$\pl\bar{X}\x_0\pl\bar{X}$, we remark that, up to a smooth half-density
section, it is a $L^1$ function for $\Re(\la)<\ndemi$,
$\la\notin \mc{R}\cup Z^1_-\cup Z^2_-$
which can be meromorphically continued to
$\cc\setminus (Z^1_-\cup Z^2_-)$
in the distribution sense (see \cite[Th. 3.2.4]{Ho} or
\cite{GZ3} for instance). This continuation makes appear some
poles in the physical sheet at $Z^2_+$ and possibly at $Z_+^1$ though
$k_1(\la)$ is holomophic at these points. The recent work
of Graham-Zworski \cite{GRZ} gives a nice description of these
poles, these are first order poles for $S(\la)$ whose residues
can be calculated explicitly. For $k\in\nn$, the residue
of $S(\la)$ at $\frac{n+k}{2}$ is the sum of
a differential operator on $\pl\bar{X}$
and of a smoothing finite-rank operator which only appears
when $\frac{n+k}{2}\in\sigma_{pp}(\Delta_g)$.
This differential operator only depends on the $k$ first
derivatives of the metric at the boundary and it is
never zero for $k$ even but can be zero for $k$ odd according to
whether the metric is even or not (it will be detailed later).

For $\la\in\cc\setminus(\mc{R}\cup \demi\zz)$,
$s(\la)$ is a polyhomogeneous conormal
distribution of order $-2\la$ associated to
$\delta_{\pl\bar{X}}$, thus $S(\la)$
is a one-step pseudodifferential operator of order $2\la-n$ on $\pl\bar{X}$.
Following Shubin's definition \cite[Def. 11.2]{SH}, $S(\la)$ is a
holomorphic family in
\[\left\{\Re(\la)<\ndemi\right\}\setminus (\mc{R}\cup Z^1_-\cup Z^2_-)\]
of zeroth order pseudodifferential operators on the compact manifold
$\pl\bar{X}$.
Therefore, $S(\la)$ is holomorphic in the same open of
$\cc$ with values in $\mc{L}(L^2(\pl\bar{X}))$.

If $h_0:=x^2g|_{T\pl\bar{X}}$, the principal symbol of $S(\la)$ is given by Joshi and
S\'a Barreto \cite{JSB}:
\beq\label{symboleprincipal}
\sigma_{0}\left (S(\la)\right )=c(\la)\sigma_{0}\left
(\Lambda^{2\la-n}\right )
\eeq
\[\Lambda:=(1+\Delta_{h_0})^\demi, \quad
c(\la):=2^{n-2\la}\frac{\Gamma(\ndemi-\la)}{\Gamma(\la-\ndemi)}\]
which leads us to set the factorization (see \cite{P1,GZ3,PP} for a similar
approach)
\[\til{S}(\la):=c(n-\la)\Lambda^{-\la+\ndemi}S(\la)\Lambda^{-\la+\ndemi}.\]
So $\til{S}(\la)$ can be expressed by
\[\til{S}(\la)=1+K(\la)\]
where $K(\la)$ is a compact operator for
$\la\in\cc\setminus (\mc{R}\cup \demi\zz)$. Notice that
the poles of $S(\la)$ on $Z^2_+$ already appear on
the principal symbol (\ref{symboleprincipal}), but
$K(\la)$ is regular on $Z_+^2$ in view of the factorization
by the Gamma factor.

Recall the functional equations satisfied by $S(\la)$ and $\til{S}(\la)$
on $\{\Re(\la)=\ndemi\}\setminus\{\ndemi\}$, (cf. \cite{GRZ})
\beq\label{eqfonct}
S^{-1}(\la)=S(n-\la)=S(\la)^* , \quad \til{S}^{-1}(\la)=\til{S}(n-\la)=\til{S}(\la)^*
\eeq
which show that $S(\la)$ is regular on the critical line $\{\Re(\la)=\ndemi\}$.

In order to use analytic Fredholm theorem to invert
$1+K(\la)$, we give the meromorphic properties of $K(\la)$,
naturally inherited from $S(\la)$, in a neighbourhood
of the physical sheet $\mc{O}_0$:

\begin{lem}\label{pseudoholo}
For all $\eps\in (0,\demi)$ and $\alpha>0$:
\beq\label{sdeltamoins1}
K(\la)\in \mc{H}ol\left (\mc{O}_\eps\setminus (Z^1_+\cup \mc{R}),
\mc{L}(L^2(\pl\bar{X}),H^{1-\alpha}(\pl\bar{X}))\right )
\eeq
At each $\la_j:=\frac{n+1}{2}+j\in Z^1_+$ with $j\in\nn_0$,
$K(\la)$ is either regular or has a first order pole.
\end{lem}
\textsl{Proof}: let $\Psi^m(\pl\bar{X},\Gamma^\demi)$
the space of pseudodifferential operator of order $m$
on $\pl\bar{X}$, acting on half densities.
As said before, the distribution
$(\beta_{\pl})_*(r^{-2\la}k_1(\la))$ is holomorphic
in $\mc{O}_\eps\setminus(Z^1_+\cup Z^2_+\cup \mc{R})$
with at most some first order poles at $Z_+^2\cup Z_+^1$.
We then kill the poles of $S(\la)$ at $Z^2_+$
by dividing by $\Gamma(\ndemi-\la)$ and according to
Shubin's definition \cite[Def. 11.2]{SH}, the operators
$(\Gamma(\ndemi-\la))^{-1}S(\la)$ form a holomorphic family
\[\frac{S(\la)}{\Gamma\left(\ndemi-\la\right)}\in
\mc{H}ol\left(\mc{O}_\eps\setminus (Z^1_+\cup\mc{R}),
\Psi^{2\Re(\la)+n+\alpha}(\pl\bar{X},\Gamma^\demi)\right), \quad
\forall\alpha>0.\]
It can easily be checked in charts by studying
its total local symbol defined by local Fourier transformation
of its Schwartz kernel,
moreover $S(\la)$ is a one-step pseudodifferential operator of
order $2\la-n$.

Therefore, the holomorphic properties of $\Lambda^{-\la+\ndemi}$
(cf. \cite[Th. 1.11]{SH}), the calculus (\ref{symboleprincipal})
and the composition properties of holomorphic pseudodifferential
operators imply that
\[K(\la)\in \mc{H}ol\left (\mc{O}_\eps\setminus (Z^1_+\cup \mc{R}),\Psi^{-1+\alpha}
(\pl\bar{X},\Gamma^\demi)\right ),\quad \forall \alpha>0.\]
At last, (\ref{sdeltamoins1}) follows immediately by using the Sobolev
continuity properties of holomorphic pseudodifferential operators
on compact manifolds.
\qed\\

We have seen that the poles of $S(\la)$ on $Z_+^2$
are somewhat artificial because they all can be captured
by the term $\Gamma(\ndemi-\la)$ of the principal symbol of $S(\la)$
and they will not be a barrier to use Fredholm theorem
in order to invert $S(\la)$, contrary to the poles on $Z^1_+$.
This follows from the fact that the scattering operator
of the model associated to this geometry
(i.e. the hyperbolic space $\hh^{n+1}$)
has some poles on $Z_+^2$ but not on $Z_+^1$.\\

\textsl{Remarks}: it is more correct to define $E(\la)$ and $S(\la)$
on sections of some conormal bundles to $\pl\bar{X}$ to keep an invariant
definition with respect to the choice of the boundary defining function $x$.
That is dropped for notations, because it does not play an
important role for what we study.

We can also notice that S\'a Barreto \cite{SB} has recently adapted the
radiation fields theory to this setting, which provides a new definition for
the operators $E(\la), S(\la)$ in terms of these radiation fields.

\subsection{Meromorphic equivalences}
Let $U\subset\cc$ an open set and $(\mc{B}_i)_{i=0,1,2}$ some Banach spaces.
Then the product $(A(\la),B(\la))\to A(\la)B(\la)$ on
\[\mc{H}ol(U,\mc{L}(\mc{B}_i,\mc{B}_j))\subset
\mc{M}er_f(U,\mc{L}(\mc{B}_i,\mc{B}_j))\subset
\mc{M}er(U,\mc{L}(\mc{B}_i,\mc{B}_j))\]
has the following stability properties
\beq\label{compomero}
\begin{array}{rcl}
\mc{M}er(U,\mc{L}(\mc{B}_i,\mc{B}_j))\x
\mc{M}er(U,\mc{L}(\mc{B}_j,\mc{B}_k))&\to& \mc{M}er(U,\mc{L}(\mc{B}_i,\mc{B}_k))\\
\mc{M}er_f(U,\mc{L}(\mc{B}_i,\mc{B}_j))\x
\mc{M}er_f(U,\mc{L}(\mc{B}_j,\mc{B}_k))&\to& \mc{M}er_f(U,\mc{L}(\mc{B}_i,\mc{B}_k))\\
\mc{H}ol(U,\mc{L}(\mc{B}_i,\mc{B}_j))\x
\mc{H}ol(U,\mc{L}(\mc{B}_j,\mc{B}_k))&\to& \mc{H}ol(U,\mc{L}(\mc{B}_i,\mc{B}_k))
\end{array}
\eeq
Let us show the following useful lemma:
\begin{lem}\label{resultat:general}
Let $\mc{B}_0\subset\mc{B}_1$ some Banach spaces and $j
:\mc{B}_0\hookrightarrow \mc{B}_1$ the continuous inclusion. Let
$U\subset \cc$ an open set, $\la_0\in U$; if $M(\la)$ is
continuous in $U\setminus\{\la_0\}$ with values in $\mc{B}_0$ and
$j\circ M(\la)\in\mc{M}er(U,\mc{B}_1)$, then $M(\la)\in
\mc{M}er(U,\mc{B}_0)$. Furthermore, if $\mc{B}_i$ are some spaces
of continuous linear maps on Banach spaces and $j\circ
M(\la)\in\mc{M}er_f(U,\mc{B}_1)$, then $M(\la)\in
\mc{M}er_f(U,\mc{B}_0)$.
\end{lem}
\textsl{Proof}: we first show that $M(\la)\in
\mc{H}ol(U\setminus\{\la_0\},\mc{B}_0)$; it is sufficient to prove
that for all $\la_1\in U\setminus\{\la_0\}$ there exists $\eps>0$
such that \beq\label{critereinteg} \int_TM(\la)d\la=0 \eeq for all
triangle $T$ included in the open disc $\{|\la-\la_1|<\eps\}$.
Note that (\ref{critereinteg}) is satisfied when we replace
$M(\la)$ by $j\circ M(\la)$. But since $j$ is continuous, the
integral of $j\circ M(\la)$ on $T$, defined as a limit of a
Riemann sum, is exactly
\[j\circ\int_TM(\la)d\la=\int_Tj\circ M(\la)d\la.\]
We then have the desired identity (\ref{critereinteg}) since
$j$ is injective. Using the same arguments and the meromorphic assumption
on $j\circ M(\la)$ we obtain the following Laurent expansion
\[j\circ M(\la)=\sum_{k=-p}^{-1} j\circ M_k(\la-\la_0)^k +j\circ H(\la)\]
\[ M_k:=\frac{1}{2\pi
i}\int_T(\la-\la_0)^{-k-1} M(\la)d\la,
\quad H(\la):=\frac{1}{2\pi i}\int_T(z-\la)^{-1} M(z)dz\]
for $\la$ near $\la_0$ and $T$ a triangle around $\la_0$. Since $j$ is
injective we have
\[M(\la)=\sum_{k=-p}^{-1}M_k(\la-\la_0)^k +H(\la)\]
and $H(\la)$ is holomorphic near $\la_0$ with values in
$\mc{B}_0$. It remains to remark that it if $j\circ M_k$ has
finite rank then it is the same for $M_k$, $j$ being injective.
\qed\\

Let us study the relations between the meromorphic properties of
 $R(\la)$, $S(\la)$ and $\til{S}(\la)$.

\begin{prop}\label{inclusion:des:poles}
Let $U\subset \{\Re(\la)<\ndemi\}$ an open set in $\cc$, the following
assertions are equivalent:\\
(1) $R(\la)$ is meromorphic in $U$.\\
(2) $S(\la)$ is meromorphic in $U$.\\
(3) $\til{S}(\la)$ is meromorphic in $U$.\\
If $U\cap Z^2_-=\varnothing$, these assertions are equivalent:\\
(4) $R(\la)$ is finite-meromorphic in $U$.\\
(5) $S(\la)$ is finite-meromorphic in $U$.\\
(6) $\til{S}(\la)$ is finite-meromorphic in $U$.\\
(7) $\til{S}(n-\la)$ is finite-meromorphic in $U$.\\
If $U\cap Z^2_-\not=\varnothing$ we just have
\[(7)\iff (6)\Rightarrow (5)\iff (4).\]
\end{prop}
\textsl{Proof}: $(2)\Rightarrow(1)$ and $(5)\Rightarrow(4)$:
let $N>\ndemi$, we first show that in the open set
\[\{\la\in \cc;  n-N< \Re(\la)< \ndemi \textrm{ and } \la,n-\la\notin
(Z^1_-\cup Z^2_-\cup \mc{R})\}\] we have the following holomorphic
identity on $\mc{L}(x^NL^2(X),x^{-N}L^2(X))$:
\beq\label{mesurespect}
R(\la)-R(n-\la)=(2\la-n)\trans E(n-\la)
S(\la)E(n-\la).
\eeq
Observe that the proof of Green's formula
obtained by Agmon \cite{A1}, Perry \cite{P} or Guillop\'e
\cite{G1} for hyperbolic quotients remains true in our framework
(see also Borthwick \cite[Prop. 4.5]{B} in our setting): for $\la,
n-\la\notin(\mc{R}\cup Z^1_-\cup Z^2_-)$, $m,m'\in X$ and $m\not=m'$
\beq\label{formulegreen}
r(\la;m,m')-r(n-\la;m,m')=
(n-2\la)\int_{\pl\bar{X}}e(\la;.,m)e(n-\la;.,m')
\eeq
which can be reformulated by
\beq\label{formulegreen2}
R(\la)-R(n-\la)=(n-2\la)\trans E(\la)E(n-\la)
\eeq
considered as
continuous operators from $\dot{C}^\infty(\bar{X},\Gamma_0^\demi)$
to $C^{-\infty}(\bar{X},\Gamma_0^\demi)$. From
(\ref{formulegreen}), (\ref{noyaudiffusion}) and
(\ref{eisensteinfct}) it is straightforward to check that
\[e(\la;y,m')=-\int_{\pl\bar{X}}s(\la;.,y)e(n-\la;.,m')\]
where $y\in\pl\bar{X}$ and $m'\in X$.
This identity can be expressed by
\beq\label{relationebla}
E(\la)=-\trans S(\la)E(n-\la)
\eeq
considered as continuous operators from
$\dot{C}^{\infty}(\bar{X},\Gamma_0^\demi)$
to $C^{\infty}(\pl\bar{X},\Gamma^\demi)$.
We deduce from (\ref{formulegreen2}) and (\ref{relationebla}) the weak identity
(\ref{mesurespect}) in the open set of $\cc$
\[\{\la\in \cc; \Re(\la)< \ndemi \textrm{ and } \la,n-\la\notin
(Z^1_-\cup Z^2_-\cup \mc{R})\}\]
Moreover, according to (\ref{propelambda}),
the fact that $S(\la)$ is holomorphic on $\mc{L}(L^2(\pl\bar{X}))$ in this open set and
(\ref{compomero}), we have the desired holomorphic identity (\ref{mesurespect}).

Let $\la_0\in \{\Re(\la)<\ndemi\}$ and $N>|\Re(\la_0)|+n$,
the identity (\ref{mesurespect}) holds near $\la_0$ with values in
$\mc{L}(x^NL^2(X),x^{-N}L^2(X))$.
On $\Re(\la)<\ndemi$, we have seen that $R(n-\la)$, $E(n-\la)$ and
$\trans E(n-\la)$ are finite-meromorphic
with only poles the points $\la_e\in\cc$ such that
$\la_e(n-\la_e)\in\sigma_{pp}(\Delta_g)$. One deduces that
(\ref{compomero}) and (\ref{mesurespect}) prove
$(2)\Rightarrow(1)$ and $(5)\Rightarrow(4)$.\\

$(1)\Rightarrow (2)$:
let $\la_0\in\{\Re(\la)<\ndemi\}$ be a pole of $R(\la)$
and $U:=B(\la_0,\eps)\subset\{\Re(\la)<\ndemi\}$ be an open disc
of $\cc$ around  $\la_0$ with radius $\eps$ taken sufficiently small to avoid
other poles of $R(\la)$.
As claimed before, $S(\la)$ is holomorphic in $U$
with values in $\mc{L}(L^2(\pl\bar{X}))$, more precisely it is a holomorphic family
of pseudodifferential operators of negative order.
In $U$, $\la_0$ is the only pole of
$x^{-\la+\ndemi}R(\la)x^{-\la+\ndemi}$ defined as an operator of
$\mc{L}(x^{2\eps} L^2(X),x^{-2\eps}L^2(X))$.
We then have the expansion in $U$
\beq\label{deverla}
x^{-\la+\ndemi}R(\la)x^{-\la+\ndemi}=\sum_{i=-p}^{-1}(\la-\la_0)^iA_i+H(\la),
\eeq
\[A_i\in \mc{L}(x^{2\eps} L^2(X),x^{-2\eps}L^2(X)),
\quad H(\la)\in\mc{H}ol(U,\mc{L}(x^{2\eps}
L^2(X),x^{-2\eps}L^2(X))).\]
Moreover the Schwartz kernels $a_k$ and $h(\la)$ of $A_k$ and $H(\la)$ can be
described by an integral on the circle $C(\la_0,\frac{\eps}{2})$
\[a_k=\frac{1}{2\pi i}\int_{C(\la_0,\frac{\eps}{2})}
(z-\la_0)^{-k-1}(xx')^{-z+\ndemi}r(z)dz,\]
\beq\label{hla}
h(\la)=\frac{1}{2\pi i}\int_{C(\la_0,\frac{\eps}{2})}
(z-\la)^{-1}(xx')^{-z+\ndemi}r(z)dz , \quad
|\la-\la_0|<\frac{\eps}{2}.
\eeq
The structure of $r(\la)$ implies
that $\beta^*(a_k)$ is the sum of a section
in $\beta^*((xx')^\ndemi C^\infty(\bar{X}\x\bar{X},\ddens))$ and of
smooth section on $(\bar{X}\x_0\bar{X})\setminus
\mc{F}$ which has a conormal singularity (not necessarily
polyhomogeneous) on $\mc{F}$ of order $-2\Re(\la_0)+n-\eps$.
Likewise $h(\la)$ is the sum of $(xx')^{-\la+\ndemi}r_0(\la)$ and a
distribution $h_1(\la)$ whose lift $\beta^*(h_1(\la))$ has the
same structure than $\beta^*(a_k)$ ($h_1(\la)$ is an integral like
(\ref{hla}) with $r_1(\la)+r_2(\la)$ instead of $r(\la)$).

Using the representation (\ref{noyaudiffusion}) of $S(\la)$
by its Schwartz kernel $s(\la)$ and the fact
that $s(\la)$ is, up to a smooth half-density section,
a $L^1$ function on $\pl\bar{X}\x\pl\bar{X}$ for $\la\in U$, we find
\[\frac{s(\la)}{2\la-n}=\sum_{k=-p}^{-1}\frac{(\la-\la_0)^k}{2\pi i}
\int_{C(\la_0,\frac{\eps}{2})}\frac{s(z)}{(2z-n)(z-\la_0)^{k+1}}dz+
\frac{1}{2\pi i}\int_{C(\la_0,\frac{\eps}{2})}
\frac{s(z)}{(2z-n)(z-\la)}dz.\]
Since $S(\la)$ is holomorphic on $\{\la\in\cc;0<|\la-\la_0|<\frac{\eps}{2}\}$
with values in $\mc{L}(L^2(\pl\bar{X}))$, we have
\[\frac{S(\la)}{2\la-n}=\sum_{k=-p}^{-1}\frac{(\la-\la_0)^k}{2\pi i}
\int_{C(\la_0,\frac{\eps}{2})}\frac{S(z)}{(2z-n)(z-\la_0)^{k+1}}dz+
\frac{1}{2\pi i}\int_{C(\la_0,\frac{\eps}{2})}
\frac{S(z)}{(2z-n)(z-\la)}dz\]
on the same set with values in $\mc{L}(L^2(\pl\bar{X}))$.
The second integral being holomorphic near $\la_0$, we conclude that
$S(\la)$ admits a finite Laurent expansion at $\la_0$.
$(1)\Rightarrow(2)$ is then proved.\\

$(4)\Rightarrow (5)$: following what we did before, it suffices
to show that if the polar part of $R(\la)$ has a finite total rank
then it is the same for $S(\la)$. Suppose (\ref{deverla})
where $A_i$ are some finite rank operators.
The Schwartz kernel of $A_i$ can be expressed by
\[a_i(x,y,x',y')=\sum_{j=1}^{r_i} \psi_{ij}(x,y)\varphi_{ij}(x',y')
\left|\frac{dxdydx'dy'}{x^{n+1}x'^{n+1}}\right|^\demi,
\quad \psi_{ij},\varphi_{ij}\in x^{-2\eps}L^2(X,dvol_g)\]
\[\dim\textrm{Vect}\{\varphi_{ij}; j=1,\dots,r_i\}=
\dim\textrm{Vect}\{\psi_{ij}; j=1,\dots,r_i\}=r_i=\rang A_i,\]
Note that elliptic regularity implies that $\psi_{ij}$ and $\varphi_{ij}$
are smooth in $X$.
Since $(\psi_{ij})_j$ are independent, one can easily see that
there exist $r_i$ points $m_1,\dots,m_{r_i}\in X$ such that the matrix
$(M_{jk})_{j,k}:=(\psi_{ij}(m_k))_{j,k}$
has rank $r_i$. Moreover
\[\phi_{ij}(x,y):=\sum_{j=1}^{r_i}\psi_{ij}(m_k)\varphi_{ij}(x,y)\in
x^\ndemi C^\infty(\bar{X})\]
since $a_i\in (xx')^\ndemi C^\infty(\bar{X}\x\bar{X}\setminus
\delta_{\pl\bar{X}},\Gamma_0^\demi)$. But $(\phi_{ij})_{j=1,\dots,r_i}$ is
a basis of $\textrm{Vect}\{\varphi_{ij}; j=1,\dots,r_i\}$, hence
\[\varphi_{ij}\in x^\ndemi C^\infty(\bar{X}),\quad j=1,\dots,r_i\]
By the same arguments, the same result holds about
$\psi_{ij}$ but with others $m_k\in X$.
The restriction of $a_i$ on $x=x'=0$ is then explicit
and $S(\la)$ can be expressed by
\[S(\la)=(2\la-n)\sum_{i=-p}^{-1}(\la-\la_0)^i\sum_{j=1}^{r_i}\psi^\sharp_{ij}
\otimes\varphi^\sharp_{ij}+H_1^\sharp(\la)\]
where $\psi^\sharp_{ij}, \varphi^\sharp_{ij}\in C^\infty(\pl\bar{X},\Gamma^\demi)$
are defined by
\[\psi^\sharp_{ij}=\left(\psi_{ij}\left|\frac{dxdy}{x^{n+1}}\right|^\demi\right)|_{x=0},\quad
\varphi^\sharp_{ij}=\left(\varphi_{ij}\left|\frac{dxdy}{x^{n+1}}\right|^\demi\right)|_{x=0}\]
and $H^\sharp_1(\la)$ is holomorphic near $\la_0$ in $\mc{L}(L^2(\pl\bar{X}))$.\\

$(2)\iff(3)$: it is sufficient to observe that
$c(n-\la)\Lambda^{-\la+\ndemi}$ and $\Lambda^{-\la+\ndemi}$ with their inverse are
meromorphic in $\mc{L}(H^p(\pl\bar{X}),H^{p-N}(\pl\bar{X}))$
for all $p\in\rr$ and $N>-\Re(\la)+\ndemi$, and conclude with
(\ref{compomero}) and Lemma \ref{resultat:general}.\\

$(6)\iff(7)$: assume that $\til{S}(\la)=1+K(\la)$ is finite-meromorphic in $U$
with $K(\la)$ compact, analytic Fredholm theorem then proves
that $\til{S}^{-1}(\la)=\til{S}(n-\la)$ is finite-meromorphic in $U$.
The reverse is identical.\\

$(6)\Rightarrow (5)$: observe that
$c(\la)\Lambda^{\la-\ndemi}$ and $\Lambda^{\la-\ndemi}$ are
holomorphic with values in $\mc{L}(L^2(\pl\bar{X}))$
and use (\ref{compomero}).\\

$(5)\Rightarrow (6)$: if $U\cap Z^2_-=\varnothing$,
$c(n-\la)\Lambda^{-\la+\ndemi}$ and $\Lambda^{-\la+\ndemi}$ are
holomorphic in $U$ with values in Sobolev spaces, so it remains
to apply (\ref{compomero}) and Lemma \ref{resultat:general}.
\qed\\

\section{Proofs of the results}

\textsl{Proof of Proposition \ref{amelioration}}: using Proposition
\ref{inclusion:des:poles} it suffices to study the meromorphic
properties of $\til{S}(\la)$. Note that (\ref{eqfonct})
implies that $\til{S}(\la)$ is unitary on $\{\Re(\la)=\ndemi\}$
(and $\la\not=\ndemi)$, so $1+K(\la)$ is invertible at one point
of $\mc{O}_\eps$. The analytic Fredholm theorem allows to prove that
\[ \til{S}^{-1}(\la)=(1+K(\la))^{-1}\]
is a well defined  finite-meromorphic family of operator in
$\mc{O_\eps}\setminus (Z^1_+\cup \mc{R})$.
Moreover since $\til{S}(n-\la)=\til{S}(\la)^{-1}$ on $\{\Re(\la)=\ndemi\}$,
the meromorphic extension of $\til{S}(\la)$ is given by
\[\til{S}(\la):=(1+K(n-\la))^{-1}\in \mc{M}er_f\left (\{\Re(\la)<\ndemi+\eps\}
\setminus (Z^1_-\cup (n-\mc{R})),
\mc{L}\left(L^2(\pl\bar{X})\right)\right ).\]
It is clear that $(n-\mc{R})\cap \{\Re(\la)<\ndemi\}$ is
the discrete spectrum, i.e. the points $\la_e$ such that
$\la_e(n-\la_e)\in \sigma_{pp}(\Delta_g)$. If $\la_e$ is
one of these `eigenvalues' which is not in $Z^1_-$,
then $n-\la_e$ is a first order pole of
finite multiplicity of $R(\la)$ and at most a first order pole
of finite multiplicity of $K(\la)$ in view of (\ref{residuvp}),
(\ref{noyaudesla}) and the fact that $S(\la)$
has at most some first order poles on $Z^2_+$.
One concludes that $\til{S}(\la)$ is finite-meromorphic
in $\{\Re(\la)< \ndemi\}\setminus Z^1_-$, hence by Proposition
\ref{inclusion:des:poles} $R(\la)$ is finite-meromorphic in
the same open set.\\

Using Proposition \ref{inclusion:des:poles}, we deduce that
$R(\la)$ is finite-meromorphic on $\mc{O}_N$ if and only if the
points of  $Z^1_-\cap \mc{O}_N$ are some poles of finite total
polar rank of $\til{S}(n-\la)$, or equivalently if and only if the
points of $Z^1_+\cap (n-\mc{O}_N)$ are some poles of finite total
polar rank of $K(\la)$. Let us show the following lemma which is a
direct consequence of the results of Graham and Zworski
\cite{GRZ}:
\begin{lem}\label{parite:et:poles}
Let $k\in\nn$, $(X,g)$ an asymptotically hyperbolic manifold and
suppose that $g$ is an even metric modulo $O(x^{2k+1})$ which we
write in the collar $U:=(0,\eps_x)\x\pl\bar{X}$
\beq\label{metriquepairemodulo} g=x^{-2}\left(dx^2+\sum_{i=0}^k
h_{2i}x^{2i}+h_{2k+1}x^{2k+1}+O(x^{2k+2})\right)
\eeq
with $(h_{2i})_{i=0,\dots,k}$ and $h_{2k+1}$ some symmetric tensors on
$\pl\bar{X}$. Then for $j=0,\dots,k+1$ the points
$\la_j:=\frac{n+1}{2}+j$ are at most some first order poles of
$S(\la)$ whose residue is given by
\[\textrm{Res}_{\la_j}S(\la)=\Pi_{\la_j}, \quad j=0,\dots,k-1\]
\beq\label{residuenlak}
\textrm{Res}_{\la_k}S(\la)=\Pi_{\la_k}-\frac{(n-\la_k)}{4}
\tra(h_0^{-1}h_{2k+1})
\eeq
\[\textrm{Res}_{\la_{k+1}}S(\la)=\Pi_{\la_{k+1}}-p_{2k+3}\]
where $\Pi_{\la_j}$ is the finite rank operator whose Schwartz kernel is
\[\pi_{\la_j}:=(2\la_j-n)\left((xx')^{-\la_j+\ndemi}
\textrm{Res}_{\la_j}r(\la)\right)
|_{\pl\bar{X}\x\pl\bar{X}},\quad j=0,\dots,k\]
$h_0^{-1}$ is the metric induced by $h_0$ on $T^*\pl\bar{X}$ and
$p_{2k+3}$ is a differential operator of order $2$ on $\pl\bar{X}$
whose principal symbol vanishes identically only if
\beq\label{conditionsymbnul}
(2-n(n-\la_k))\tra(h_0^{-1}h_{2k+1})=0
\eeq
\end{lem}
\textsl{Proof}: first note that for $m\in \pl\bar{X}$
the tensors $h_i(m)$ will be considered as symmetric matrices
in $\rr^n$ via the Euclidean scalar product of the chart,
 $\tra(h_0^{-1}h_{2k+1})$ can be understood in that way or
by the trace of the linear operator associated to $h_{2k+1}$ via the
scalar product $h_0$ on $T\pl\bar{X}$.
Now we use the construction of the Poisson operator
according to \cite{GRZ}. The first step is to construct,
for a given $f_0\in C^\infty(\pl\bar{X})$, a solution $F_\infty\in C^\infty(\bar{X})$ of
\beq\label{solutionapprochee}
(\Delta_g-\la(n-\la))x^{n-\la}F_\infty=O(x^\infty), \quad F_\infty|_{x=0}=f_0
\eeq
it is then clear that a solution in the collar $U_x$ is sufficient (we
can always multiply it by a smooth cut-off function with support near
$\pl\bar{X}$ and which is equal to $1$ in a neighbourhood of $\pl\bar{X}$).
Let $M_t=\{x=t\}$ and $h(t)$ (resp. $h^{-1}(t)$)
the metric induced by $x^2g$ on $TM_t$ (resp. $T^*M_t$).
In $U_x$, we first have the identity
\[(\Delta_g-\la(n-\la))x^{n-\la}=x^{n-\la}\mc{D}_\la\]
\[\mc{D}_\la:=-x^2\pl_x^2+\left(2\la-n-1-\frac{x}{2}\tra(h^{-1}(x)\pl_xh(x))\right)
x\pl_x-\frac{(n-\la)x}{2}\tra(h^{-1}(x)\pl_xh(x))+x^2\Delta_{h(x)}\]
and for $f\in C^\infty(\pl\bar{X})$ and $j\in\nn_0$
\beq\label{relationdla}
\mc{D}_\la(fx^j)=j(2\la-n-j)fx^{j}+x^{j}G(\la-j)f,
\eeq
\[(G(z)f)(x,y):=x^2\Delta_{h(x)}f(y)-
\frac{(n-z)x}{2}\tra(h^{-1}(x)\pl_xh(x))f(y).\]
Suppose now that $F\in C^\infty(\bar{X})$ is a function
such that $\mc{D}_\la F=O(x^j)$ (with $j\geq 1$), then
since $G(z)f=O(x)$ for all $f\in C^\infty(\pl\bar{X})$,
(\ref{relationdla}) ensures that
\[\mc{D}_\la F -\mc{D}_\la \left(x^j
\frac{(x^{-j}\mc{D}_\la(F))|_{x=0}}{j(2\la-n-j)}\right)=O(x^{j+1}).\]
Let $f_0\in C^\infty(\pl\bar{X})$ fixed; since (\ref{relationdla})
implies that $\mc{D}_\la f_0=O(x)$, the previous remark allows us to construct
the functions $F_j\in C^\infty(\bar{X})$
(for $j\geq 0$) and $f_j\in C^\infty(\pl\bar{X})$ (for $j\geq 1$) by the
induction formula
\beq\label{constructionfj}
F_0=f_0, \quad f_j=\frac{-(x^{-j}\mc{D}_\la(F_{j-1}))|_{x=0}}{j(2\la-n-j)},\quad
F_j=F_{j-1}+f_jx^j, \quad j\geq 1.
\eeq
By construction, we obtain
\[\mc{D}_\la(F_{j-1})=O(x^j)\]
and according to Borel lemma, there exists a function
$F_\infty\in C^\infty(\bar{X})$ whose Taylor coefficients
at $x=0$ are $f_j$, which gives that $F_\infty$ is a
solution for (\ref{solutionapprochee}).
In other respects, we can express $f_j$ by
\[f_j=p_{j,\la}f_0\]
where $p_{j,\la}$ is a differential operator of order $\leq 2[j/2]$
on $\pl\bar{X}$.

Recall that Proposition 3.5 of \cite{GRZ} proves that
the residues of $S(\la)$ at $(\la_l)_{l=0,\dots,k+1}$ are
\beq\label{propzworski}
\textrm{Res}_{\la_l}S(\la)=\Pi_{\la_l}-p_{2l+1},
\quad p_{2l+1}:=\textrm{Res}_{\la_l}(p_{2l+1,\la}).
\eeq
Consequently, it remains to calculate $(p_{2l+1})_{l=0,\dots,k+1}$.\\

We will denote by $D^j$ the set of differential operators of order
$j$ on $\pl\bar{X}$ and for notational simplicity $D^j$ also means
all differential operator of order $j$ on $\pl\bar{X}$ that we do
not need to know explicitly. Let us now set
\[K:=\tra(h_0^{-1}h_{2k+1}).\]
In the Taylor expansion of $G(z)$ at $x=0$, we use the assumption
(\ref{metriquepairemodulo}) and group the even
powers of $x$ together in $G_2(z)$ and the odd powers of $x$
together in $G_1(z)$ to obtain $G(z)=G_1(z)+G_2(z)$ with
\begin{eqnarray}
\label{g1}
G_1(z)&=&-x^{2k+1}\frac{(n-z)(2k+1)}{2}K+x^{2k+3}Q+O(x^{2k+5}),\\
\label{g2} G_2(z)&=&x^2\Delta_{h_0}+x^2D^0+O(x^4),
\end{eqnarray}
\[ Q\in D^2, \quad
\sigma_0(Q)(\xi)=-\cjg h_0^{-1}h_{2k+1}h_0^{-1}\xi,\xi\cjd,\]
where $\sigma_0(Q)$ is the principal symbol (of order $2$) of $Q$.
Hence, a first application is that for all $f\in C^\infty(\pl\bar{X})$
\beq\label{impairemodulo}
\mc{D}_\la(x^{2j}f) \textrm{ is even modulo } O(x^{2k+2j+1}).
\eeq

We then show by induction that $F_j$ is even in $x$
for $j=0,\dots,2k$. $F_0$ is even, suppose now
that $F_{j-1}$ is even for a fixed $j\leq 2k-1$. If $j$ is even,
(\ref{constructionfj}) clearly implies that
$F_j$ is even. On the other hand, $F_{j-1}$ being even,
(\ref{impairemodulo}) shows that $\mc{D}_\la(F_{j-1})$
is even modulo $O(x^{2k+1})$. So if $j$ is odd,
$x^{-j}\mc{D}_\la(F_{j-1})$ is odd modulo $O(x^{2k+1-j})$
and it vanishes at $x=0$ since $2k+1-j\geq 2$ by assumption on $j$.
From the definitions of $f_j$ and $F_j$ in
(\ref{constructionfj}) we finally obtain that $f_j=0$ and that
$F_j=F_{j-1}$ is even.
As a conclusion $p_{2l+1,\la}=p_{2l+1}=0$ for $l=0,\dots,k-1$.\\

Now, to construct $f_{2k+1}$, remark
(\ref{impairemodulo}) shows that the coefficient of order $2k+1$
of $\mc{D}_\la F_{2k}$ is exactly the coefficient of order $2k+1$
of $\mc{D}_\la f_0$, namely
\[-\frac{(n-\la)(2k+1)}{2}Kx^{2k+1}f_0.\]
One concludes that
$p_{2k+1,\la}=\frac{1}{4}(n-\la)(\la-\la_k)^{-1}K$, hence its
residue at $\la=\frac{n+1}{2}+k$ is
\[p_{2k+1}=\frac{(n-\la_k)}{4}K \]
and (\ref{residuenlak}) is obtained using (\ref{propzworski}).\\

For the term $p_{2k+3,\la}$, we shall only study its principal
symbol.
To obtain $f_{2k+3}$, we need to evaluate the coefficient before
$x^{2k+3}$ in $\mc{D}_\la(\sum_{i=0}^{2k+2}x^if_i)$. But since
$f_{2l+1}=0$ for $l<k$ we can use (\ref{impairemodulo}) to check
that the only terms having a non zero coefficient before
$x^{2k+3}$ in $\mc{D}_\la(F_{2k+2})$ come from $\mc{D}_\la f_i$
with $i\in\{0,2,2k+1\}$.

Consider now the three cases. According to (\ref{g1}), the term
of order $x^{2k+3}$ in $\mc{D}_\la f_0$ is
\[Qf_0x^{2k+3}.\]
The one in $\mc{D}_\la f_2$ is
\[-\demi(n-\la+2)(2k+1)K f_2
x^{2k+3}=\frac{(n-\la+2)(2k+1)}{4(2\la-n-2)}K
(\Delta_{h_0}+D^0)f_0 x^{2k+3}.\]
Finally for $\mc{D}_\la f_{2k+1}$, the term of order $x^{2k+3}$ comes from
$x^{2k+1}G_2(\la-2k-1)f_{2k+1}$, it is
\[(\Delta_{h_0}+D^0)f_{2k+1}x^{2k+3}=
\frac{(n-\la)}{4(\la-\la_k)}(\Delta_{h_0}+D^0)K f_0.\]
As before, let us set $\la_j:=\frac{n+1}{2}+j$ for $j\in\nn_0$.
We then deduce that
\[f_{2k+3}=\frac{-1}{2(2k+3)(\la-\la_{k+1})}\left(Q+
\frac{(n-\la+2)(2k+1)}{4(2\la-n-2)}K\Delta_{h_0}+
\frac{(n-\la)}{4(\la-\la_k)}\Delta_{h_0}K+D^0\right)f_0.\]
Now taking the residue at $\la_{k+1}$ we find
\[p_{2k+3}=-\frac{1}{2(2k+3)}\left(Q+\frac{(n-\la_k)}{2}K\Delta_{h_0}+D^1\right)\]
which is a differential operator with principal symbol
\beq\label{symbprincipal}
\sigma_0(p_{2k+3})(\xi)=-\frac{1}{2(2k+3)}\left\cjg
-h_0^{-1}h_{2k+1}h_0^{-1}\xi+\frac{(n-\la_k)}{2}K h_0^{-1}\xi,\xi\right\cjd.
\eeq
If $\sigma_0(p_{2k+3})=0$ we have
\[h_0^{-1}h_{2k+1}-\frac{(n-\la_k)}{2}K=0,\]
so it remains to take the trace of this identity and find (\ref{conditionsymbnul}).
\qed\\

If $g$ is even  modulo $O(x^{2k+1})$,
the residue of $S(\la)$ at $\la_j$ has finite rank for $j=0,\dots,k-1$
according to Lemma \ref{parite:et:poles}, so it is the very same thing for
$\til{S}(\la)$ and Proposition \ref{inclusion:des:poles} then proves that
$R(\la)$ is finite-meromorphic near the points $(\frac{n-1}{2}-j)_{j=0,\dots,k-1}$.

Conversely, if $R(\la)$ is finite-meromorphic near these points (with $k\geq 2$),
Proposition \ref{inclusion:des:poles} tells us that $\til{S}(\la)$ is finite-meromorphic
near the points $(\la_j)_{j=0,\dots,k-1}$. Assume that $g$ is even modulo
$O(x^{2l+1})$ with $l\leq k-2$, the residues of
$\til{S}(\la)$ at $\la_l$ and $\la_{l+1}$ must be some finite rank operators.
Following  Lemma \ref{parite:et:poles}, this implies that $(n-\la_l)\tra(h_0^{-1}h_{2l+1})=0$ and
$(2-n(n-\la_l))\tra(h_0^{-1}h_{2l+1})=0$, thus $\tra(h_0^{-1}h_{2l+1})=0$.
Taking the expression of the principal symbol of $p_{2l+3}$
in (\ref{symbprincipal}) yields $h_{2l+1}=0$, so $g$
must be even modulo $O(x^{2l+3})$. Since $g$ is always even modulo $O(x)$, an easy induction
proves that $g$ is even modulo $O(x^{2k-1})$ and the proof of Proposition \ref{amelioration} is achieved.
\qed\\

\textsl{Remark}: notice that the same kind of arguments (or using
Joshi-S\'a Barreto \cite[Th. 1.2]{JSB} formula) prove
that the set of residues of $S(\la)$ on $(\frac{n+k}{2})_{k\in\nn}$
determines the Taylor expansion of the metric at the boundary.\\

Observe that a choice of even metric modulo $O(x^{2k+1})$ satisfying
\[\textrm{meas}\left\{\tra(h_0^{-1}h_{2k+1})=0\right\}=0,
\quad \la_k\not=n\]
with the notation of (\ref{metriquepairemodulo}) gives a
residue of $S(\la)$ (and $\til{S}(\la)$) at $\la_k$ which is injective
modulo the projection on the $L^2$-eigenspace.
Let us show that $\til{S}^{-1}(\la)$ must have
an essential singularity at $\la_k$ in this case.

\begin{lem}\label{compact:injectif}
Let $\mc{B}$ a Banach space of infinite dimension,
$\la_0\in\cc$ and $U$ a neighbourhood of $\la_0$.
Let $M(\la)\in\mc{H}ol(U\setminus\{\la_0\},\mc{L}(\mc{B}))$
a meromorphic family  of bounded operators in $U$ which satisfies
\beq\label{hypothesek1}
M(\la)=1+\frac{K_{-1}}{\la-\la_0}+K(\la), \quad K(\la)\in
\mc{H}ol(U,\mc{L}(\mc{B})),
\eeq
where $K_{-1}$ and $K(\la)$ are compact and
\[\dim\ker K_{-1}<\infty.\]
If there exists $z\in U$ such that $M(z)$ is invertible,
then $M(\la)$ is invertible for almost every $\la\in U$ with inverse
$M^{-1}(\la)$ finite-meromorphic in $U\setminus\{\la_0\}$
and $\la_0$ is an essential singularity of $M^{-1}(\la)$.
\end{lem}
\textsl{Proof}: to simplify, we take $\la_0=0$.
The fact that $M(\la)$ is invertible almost everywhere in $U$
with inverse finite-meromorphic in $U\setminus\{0\}$ is a consequence
of analytic Fredholm theorem.
Assume now that  $M^{-1}(\la)$ has a finite Laurent expansion at $0$
\[M^{-1}(\la)=\sum_{i=-p}^\infty N_i\la^i , \quad p\geq 0.\]
We can take the Laurent expansion of $M(\la)$ at $0$
\[M(\la)=\sum_{i=-1}^\infty M_i\la^i=K_{-1}\la^{-1}+1+K_0+\sum_{i=1}^\infty
K_i\la^i,\]
where $K_i$ is compact, and make the product
\[M(\la)M^{-1}(\la)=\sum_{i=-p-1}^\infty\la^i\sum_{j+k=i}M_jN_k=1,\]
which leads to the system 
\beq\label{egaliteserie}
\sum_{j=-1}^{i+p}M_jN_{i-j}=\delta_{i0} , \quad i\geq-p-1. 
\eeq 
Let us show by induction that $N_i$ has finite rank for $i\leq 0$.
Taking equation (\ref{egaliteserie}) with $i=-p-1$ yields
$K_{-1}N_{-p}=0$, and by assumption on $K_{-1}$ we find that
$N_{-p}$ has finite rank. Let $I\leq -1$, suppose now that $N_i$
has finite rank for all $i\leq I$ and let us prove that $N_{I+1}$
has finite rank. For $i=I$, equation (\ref{egaliteserie}) implies
that
\[K_{-1}N_{I+1}=-\sum_{j=0}^{I+p}M_jN_{I-j}\]
has finite rank by induction assumptions (since $I-j\leq I$).
Let $r$ an integer such that
\[r>\dim\ker K_{-1}+\dim\textrm{Im}K_{-1}N_{I+1}\]
If $N_{I+1}$ has infinite rank, there exists a family of independent
vectors $(\varphi_i)_{i=1,\dots,r}$ in $\textrm{Im}(N_{I+1})$ and the
restriction of $K_{-1}$ on
$E:=\bigoplus_{i=1}^r \cc\varphi_i$ is a linear map
on a vector space of finite dimension satisfying
\[\dim\ker K_{-1}|_E+\rang K_{-1}|_E<\dim E=r\]
which is not possible. We deduce that $N_{I+1}$ has finite rank.

Take (\ref{egaliteserie}) with $i=0$: $K_{-1}$ and $(N_i)_{i\leq 0}$
being compact, it is clear that
\[1=K_{-1}N_1+\sum_{j=0}^{p}M_jN_{-j}\]
is compact, which is not possible.
\qed\\

\begin{cor}\label{naturelak}
Let $k\in\nn_0$ and $(X,g)$ an asymptotically hyperbolic manifold
and $g$ even modulo $O(x^{2k+1})$ that we write as in
(\ref{metriquepairemodulo}). If $k\not=\frac{n-1}{2}$ and
\beq\label{mesure}
\textrm{meas}\{\tra(h_0^{-1}h_{2k+1})=0\}=0,
\eeq
then $n-\la_k=\frac{n-1}{2}-k$ is an essential singularity of $R(\la)$.
\end{cor}
\textsl{Proof}: it suffices to combine Lemma
\ref{parite:et:poles} and Lemma \ref{compact:injectif} and remark
that the multiplication by a smooth function
on $L^2(\pl\bar{X})$ is injective if the measure of its zeros
vanishes. To show that the residue of $\til{S}(\la)$ at $\la_k$ has a
kernel of finite dimension, we use that the sum
of a bounded injective operator and a finite rank operator has
a finite dimensional kernel.
\qed\\

For $k=0$, $m(x,g):=(2n)^{-1}\tra(h_0^{-1}h_1)$
is exactly the mean curvature of $\pl\bar{X}$
in $(\bar{X},x^2g)$, it is a smooth function on $\pl\bar{X}$ which
depends on $x$. However it can be defined invariantly with respect to $x$
as a smooth section $m(g)$ of the conormal bundle $|N^*\pl\bar{X}|$.
In other words, a new choice of boundary defining function $t=e^\omega x$
(with $\omega\in C^\infty(\bar{X})$) gives $m(t,g)=e^{-\omega_0}m(x,g)$
where $\omega_0:=\omega|_{\pl\bar{X}}$.
If the mean curvature is almost everywhere non zero and $n\not=1$,
the corollary claims that $\frac{n-1}{2}$ is an essential singularity of $R(\la)$.

Recall that $\mc{M}_{ah}(X)$ is the space of asymptotically hyperbolic
metrics on $X$. If $x_0$ is a fixed boundary defining function, the map
\[\left\{\begin{array}{rcl}
\mc{M}_{ah}(X)&\to& \{G\in C^\infty(\bar{X},S_+^2(T^*\bar{X}));
|dx_0|_{G}=1 \textrm{ on }\pl\bar{X}\}\\
g&\to & x_0^2g
\end{array}\right.\]
is bijective and we shall identify these two spaces.
$\mc{M}_{ah}(X)$ inherits its $C^\infty$ topology from
$C^\infty(\bar{X},T^*\bar{X}\otimes T^*\bar{X})$ which is defined
as usual by semi-norms $(N_i)_{i\in\nn}$ measuring the derivatives
of the tensors ($\bar{X}$ is compact). It is not difficult to see
that the mean curvature $m(.)$ is continuous from $\mc{M}_{ah}(X)$
to
$C^\infty(\pl\bar{X},|N^*\pl\bar{X}|)$.\\

\textsl{Proof of Theorem \ref{genericitean}}: according to Thom
theorem, the set of Morse functions on $\pl\bar{X}$ is open and
dense in $C^\infty(\pl\bar{X})$ and it is exactly the same for the
Morse sections of $C^\infty(\pl\bar{X}, |N^*\pl\bar{X}|)$. Let us
denote by $V$ this subset of $C^\infty(\pl\bar{X},
|N^*\pl\bar{X}|)$. If $s\in V$ then $\textrm{meas}(s^{-1}(0))=0$,
hence $m^{-1}(V)$ is an open set of $\mc{M}_{ah}(X)$ contained in
the set of metrics in $\mc{M}_{ah}(X)$ for which $R(\la)$ has an
essential singularity at $\frac{n-1}{2}$. It remains to show that
$m^{-1}(V)$ is dense in $\mc{M}_{ah}(X)$.

We first check that $m(.)$ is a surjective map. Let
$g_0\in\mc{M}_{ah}(X)$ and take its model form
\[g_0=x^{-2}(dx^2+h_0+h_1x+O(x^2)),\]
where $x\in Z(\pl\bar{X})$. Let
\beq\label{gepsvarphi}
g_{\eps,\varphi}:=g_0+x^{-1}h_0\varphi\chi(\eps^{-1}x),
\eeq 
where $\eps>0$, $\varphi\in C^\infty(\pl\bar{X})$ and
$\chi\in C_0^\infty([0,2])$ such that $\chi(t)=1$ if $t\leq \demi$
and $\chi(t)=0$ if $t\geq 1$. $g_{\eps,\varphi}$ is an
asymptotically hyperbolic metric if $\eps$ is taken sufficiently
small (but depending on $\sup_{\pl\bar{X}}|\varphi|$).
Nevertheless $\eps$ can be chosen independent with respect to
$\varphi$ if $|\varphi|\leq 1$, we will denote it $\eps_0$. It is
straightforward to see that
\[m(g_{\eps,\varphi})=m(g_0)+\varphi|dx|,\]
hence each section $f|dx|$ of
$C^\infty(\pl\bar{X},|N^*\pl\bar{X}|)$ can be written
$m(g_{\eps,\varphi})$ by taking $\varphi:=f-m(g_0)|dx|^{-1}$
(with notation (\ref{gepsvarphi})).

Let $g_0\in \mc{M}_{ah}(X)$, $\psi_0:=m(g_0)$, $\eps_0$ defined as
before and $B(g_0):=\cap_{i\in I}B_i(g_0,r_i)$ a finite
intersection in $\mc{M}_{ah}(X)$ of `open balls' around $g_0$ with
radius $r_i$ for the semi-norms $N_i$. Let $I_0$ be the largest
number of derivatives of the metric measured by the semi-norms
$(N_i)_{i\in I}$. We set $W(\psi_0)$
a finite intersection of `open balls' around
$\psi_0$ for some semi-norms of
$C^\infty(\pl\bar{X},|N^*\pl\bar{X}|)$ which control the $I_0$
first derivatives of the section on $\pl\bar{X}$. The radius of
these balls can be chosen sufficiently small (depending on
$(r_i)_{i\in I}$) such that for all section $s\in W(\psi_0)$ the
tensor $g_{\eps_0,\varphi}$ defined in (\ref{gepsvarphi}) with
$\varphi:=(s-\psi_0)|dx|^{-1}$ lies in $B(g_0)$, in other words
the function $\varphi\to g_{\eps_0,\varphi}$ is continuous in a neighbourhood
of $0$. Note that $g_{\eps_0,\varphi}$ is a metric
because we are in the case where
$\sup_{\pl\bar{X}}|\varphi|$ approaches $0$ and we can suppose
$|\varphi|\leq 1$. Since $V$ is dense in
$C^\infty(\pl\bar{X},|N^*\pl\bar{X}|)$, one can find a Morse
section $s$ in the neibourhood $W(\psi_0)$ of $\psi_0$
and since $m(g_{\eps_0,\varphi})=s$, there exists
a metric in $m^{-1}(V)\cap B(g_0)$. We conclude that $m^{-1}(V)$
is dense in $\mc{M}_{ah}(X)$ and the proof is achieved. Concerning
$n-\la_k=\frac{n-1}{2}-k$ with $k>0$, we can show the same result
by using the residue calculus of $S(\la)$ at $\la_k$ in Lemma
\ref{parite:et:poles} for even metrics modulo $O(x^{2k+1})$.
\qed\\

\textsl{Remark}: assume that $\pl\bar{X}$ is connected and that
there exists an analytic neighbourhood of $\pl\bar{X}$ in
$\bar{X}$, hence we can take a boundary defining function $x$
which is analytic near $\pl\bar{X}$. If $x^2g$ is analytic for
$g\in\mc{M}_{ah}(X)$, then $m(g)$ is analytic on $\pl\bar{X}$ and
\[\textrm{meas}(\{m(g)=0\})>0\iff m(g)=0.\]
If $n\geq 2$ and $m(g)$ is not identically zero,
$n-\la_0=\frac{n-1}{2}$ is an essential singularity of $R(\la)$
according to Corollary \ref{naturelak}. On the other hand if
$m(g)=0$, the arguments preceding Lemma \ref{parite:et:poles}
prove that $R(\la)$ is finite-meromorphic near $n-\la_0$.
Consequently, we obtain that
\[R(\la) \textrm{ is meromorphic near } \frac{n-1}{2}\iff \pl\bar{X}
\textrm{ is minimal in } (\bar{X},x^2g),\] 
where the minimality
of $\pl\bar{X}$ does not depend on the choice of $x$.\\

\textsl{A particular case}: suppose now that $n=1$,
$g\in\mc{M}_{ah}(X)$ analytic near $\pl\bar{X}$ with $\pl\bar{X}$
connected. $S(\la)$ is holomorphic at $\la_0=1$ according to Lemma
\ref{parite:et:poles} and the fact that $0\notin
\sigma_{pp}(\Delta_g)$. Furthermore, there exists an analytic
boundary defining function $x$ such that
\beq\label{develop}
x^2g=dx^2+\sum_{i=0}^\infty h_ix^i, \quad h_i\in
C^\infty(\pl\bar{X}, S^2(T^*\pl\bar{X})) \eeq
near $\pl\bar{X}$.
If $g$ is even modulo $O(x^{2k+1})$ with $k\geq 1$, the previous
remark and Lemma \ref{parite:et:poles} imply that $R(\la)$ is
meromorphic on $\cc\setminus \cup_{i>k}\{n-\la_i\}$ if and only if
$h_{2k+1}=0$. If we consider the space of even metric modulo
$O(x^3)$ (i.e. such that $\pl\bar{X}$ is a geodesic of
$(\bar{X},x^2g)$) we obtain by induction that $R(\la)$ is
meromorphic if and only if $g$ is even and in that case the poles
have finite multiplicity.

\section{Examples with accumulation of resonances}

Before giving the example of Proposition \ref{singessentielles},
it is useful to remark that we can easily find
an asymptotically hyperbolic manifold
such that the term (\ref{residuenlak}) is a constant (not $0$).
In that case, Lemma \ref{parite:et:poles} implies that
$S(\la)$ has a first order pole of infinite multiplicity at
$\la_k=\frac{n+1}{2}+k$ and the residue is
a constant. Therefore, after renormalization, $S(\la)$
can be expressed near $\la_k$ by
\[\til{S}(\la):=1+c_k\frac{\Lambda^{-1-2k}}{\la-\la_k}+H(\la)\]
with $c_k\not=0$ and $H(\la)$ holomorphic compact.
If $(\phi_j)_{j\in\nn}$ is an orthonormal basis of eigenfunctions
of $\Lambda$ on $L^2(\pl\bar{X})$ and $(\alpha_j)_{j\in\nn}$ the associated
eigenvalues, $H(\la)\phi_j$ converge strongly to $0$ when $j\to\infty$,
thus for a large fixed $j$, $H(\la)\phi_j$ becomes insignificant
in the expression of $\til{S}(\la)\phi_j$ when $\la$ is close to $\la_k$.
If $H(\la)$ was null, we would have the following formula for the
inverse of $\til{S}(\la)$
\[\til{S}(\la)^{-1}=\sum_{j\in\nn}\frac{\la-\la_k}
{\la-\la_k+c_k\alpha_j^{-1-2k}}\phi_j\cjg.,\phi_j\cjd.\]
In other words, $(1+c_k\frac{\Lambda^{-1-2k}}{\la-\la_k})^{-1}$
has a sequence of poles $z_j=\la_k-c_k\alpha_j^{-1-2k}$ which converges to $\la_k$.
As a conclusion, the functional equation (\ref{eqfonct}) could be used
to argue that $\til{S}(\la)$ has a sequence of poles which converge to $n-\la_k$.
The key to control the `little' perturbation $H(\la)$ will be Rouch\'e's theorem.\\

\textsl{Proof of Proposition \ref{singessentielles}}: let
$k\in\nn_0$ such that $2k\not=n-1$ and the collar $U:=(0,2)\x S^n$
which carries the metric
\[g:=x^{-2}(dx^2+d(x)h_0),\]
\[d(x):=\left (1-\frac{x^2}{4}\right )^2+\chi(x)x^{2k+1},\]
where $h_0:=g_{S^n}$ is the canonical metric on the n-dimensional sphere
$S^n$, and $\chi\in C^\infty ([0,2])$ a non negative
function such that $\chi(x)=1$
for $x\in [0,\demi]$ and $\chi(x)=0$ for $x\in [1,2]$.
Set $B_{n+1}:=\{m\in\rr^{n+1};|m|<1\}$ and the
diffeomorphism
\begin{eqnarray*}
 S^n\x(0,2)&\underset{\psi^{-1}}{\to}& B_{n+1}\setminus
 \{0\}\\
(\omega,x)&\to& \frac{2-x}{2+x}\omega.
\end{eqnarray*}
It is straightforward to see that $\psi^*g$ can be extended
smoothly on $B_{n+1}$ (it is actually the hyperbolic metric
$4|dm|^2(1-|m|^2)^{2}$ in $\{|m|\leq \frac{1}{3}\}$).
Set $(X,G)$ the obtained asymptotically hyperbolic
manifold, which clearly satisfies the assumptions (\ref{metriquepairemodulo}):
we get $\hh^{n+1}$ perturbed by a $O(x^{2k+1})$ near the
boundary $S^n$.

We then obtain
\[\tra(h_0^{-1}h_{2k+1})=n\]
and we have the following expression for the Laplacian in the collar $U$
\[\Delta_g=-x^2\pl_x^2+(n-1)x\pl_x+\frac{x^2}{d(x)}\Delta_{h_0}-
\ndemi\frac{d'(x)}{d(x)}x^2\pl_x.\]
Lemma \ref{parite:et:poles} implies that the scattering operator
$S(\la)$ associated to $\Delta_G$ has a first order pole at $\la_k$ with residue
\beq\label{residusla}
\textrm{Res}_{\la_k}S(\la)=-\frac{n(n-\la_k)}{4}+\Pi_{\la_k},
\eeq
where $\Pi_{\la_k}$ is a finite rank operator on $L^2(S^n)$
whose structure can be detailed
\[\Pi_{\la_k}=\sum_{l=1}^m \varphi_l\otimes\varphi_l, \quad \varphi_l\in
C^\infty(S^n)\subset L^2(S^n).\]

Let $(v_j)_{j\in\nn}$ the eigenvalues of $\Delta_{h_0}$
(repeated with multiplicity) and $(\phi_j)_{j\in\nn}$
the associated orthonormal basis of eigenfunctions.
Since the metric $G$ is radial on $B_{n+1}$,
the Laplace operator $\Delta_G$ can be decomposed into a direct sum
\[\Delta_G\simeq \bigoplus_{j\in\nn} P_j,
\quad P_j:=-x^2\pl_{x}^2+(n-1)x\pl_x+\frac{x^2}{d(x)}v_j-\ndemi\frac{d'(x)}{d(x)}
x^2\pl_{x}\]
on
\[L^2(X,dvol_G)\simeq l^2\left (\nn,L^2\left
((0,2],\frac{d(x)^\ndemi}{x^{n+1}} dx\right )\right)\]
with singular Dirichlet condition at $x=2$.
We deduce that the resolvent and the scattering operator
can also be decomposed into a direct sum
\[R(\la)=\bigoplus_{j\in\nn} (P_j-\la(n-\la))^{-1},
\quad S(\la)=\bigoplus_{j\in\nn} S_j(\la).\]
On the other hand, recall that the expression
of the principal symbol of $S(\la)$
given in (\ref{symboleprincipal}) allows to factorize
\[c(n-\la)\Lambda^{-\la+\ndemi}S(\la)\Lambda^{-\la+\ndemi}=1+K(\la),\]
where $K(\la)$ is a meromorphic family of compact operators on $L^2(S^n)$.
Given the expression of the residue (\ref{residusla}),
we obtain for $\la$ in a neighbourhood $V_k$ of $\la_k$
\[1+K(\la)=1+\frac{K_{\la_k}}{\la-\la_k}+ H(\la),\]
\[K_{\la_k}:=c\left(n-\la_k\right)\left(-\frac{n(n-\la_k)}{4}\Lambda^{-1-2k}
+\Lambda^{-\demi-k}\Pi_{\la_k}\Lambda^{-\demi-k}\right),\]
\beq\label{prophla}
H(\la)\in\mc{H}ol(V_k,\mc{L}(L^2(S^n),H^{1-\eps}(S^n))), \quad \forall \eps>0,
\eeq
using the decomposition of $S(\la)$ on the orthonormal basis $(\phi_j)_{j\in\nn}$
of $L^2(S^n)$, we have
\[K(\la)\phi_j=K_j(\la)\phi_j, \quad K_j(\la):=\cjg K(\la)\phi_j,\phi_j\cjd,\]
\[H(\la)\phi_j=H_j(\la)\phi_j, \quad H_j(\la):=\cjg H(\la)\phi_j,\phi_j\cjd,\]
with $H_j(\la)$ holomorphic on $V_k$ and satisfying on $V_k$
\begin{eqnarray*}
|H_j(\la)|&\leq &|\cjg \Lambda^{-\demi}\Lambda^{\demi}H(\la)\phi_j,\phi_j\cjd|\\
&\leq& ||\Lambda^{\demi}H(\la)\phi_j||\alpha_j^\demi\leq C\alpha_j^\demi,
\end{eqnarray*}
where $\alpha_j:=(1+v_j)^{-\demi}$ converge to $0$ when $j\to\infty$
and $C>0$. Observe
that (\ref{prophla}) and the bound of $||\Lambda^{\demi}H(\la)||$ are
direct consequences of (\ref{sdeltamoins1}).

In the same way, $\cc\phi_j$ is globally fixed by $\Pi_{\la_k}$
\[\Pi_{\la_k}\phi_j=\beta_j\phi_j,\quad \beta_j:=\cjg\Pi_{\la_k}\phi_j,\phi_j\cjd=
\sum_{l=1}^m\cjg\varphi_l,\phi_j\cjd\cjg\phi_j,\bar{\varphi_l}\cjd\underset{j\to\infty}\to
0.\]
We deduce that
\[1+K_j(\la)=1+\frac{m_k}{\la-\la_k}\alpha_j^{1+2k}+H_j(\la),\]
\[ m_k:=c\left(n-\la_k\right)\left[\frac{n(n-\la_k)}{4}+\beta_{j}\right].\]
Remark that if $k>\frac{n-1}{2}$, then $\Pi_{\la_k}=0$
since $\la_k(n-\la_k)=\frac{n^2}{4}-(\demi+k)^2\notin \sigma_{pp}(\Delta_G)$.
We check that if $\beta_j=0$ then $m_k\not=0$, thus there exists
$J\in\nn$ such that for $j\geq J$ we have $m_k\not=0$ since $\beta_j\to 0$
when $j\to\infty$.
Consequently, we obtain an explicit expression to inverse $S_j(\la)$ for $j\geq J$
and $\la\in V_k$
\[S_j(\la)^{-1}=c(n-\la)\alpha_j^{2\la-n}
\frac{\la-\la_k}{(\la-\la_k)(1+H_j(\la))+m_k\alpha_j^{1+2k}}\]
if the denominator is not $0$.

Let us choose $\eps>0$ such that the disc with
centre $\la_k$ and radius $\eps$
is included in $V_k$ and set $z_j:=\la_k-m_k\alpha_j^{1+2k}$.
There exists an integer $J_0\geq J$ such that
the circle $C(z_j,\frac{\eps}{2})$
with centre $z_j$ and radius $\frac{\eps}{2}$ is included in $V_k$
since $z_j\to\la_k$ when $j\to\infty$. Since $J_0$ can be chosen as large as we want,
it is not restrictive to suppose that $|m_k|\alpha_j^{1+2k}\leq \eps$, so set
\[\eps_j=\frac{|z_j-\la_k|}{2}=\frac{|m_k|\alpha_j^{1+2k}}{2}\]
and both following holomorphic functions in $V_k$
\[f_j(\la):=\la-z_j=\la-\la_k+m_k\alpha_j^{1+2k},\quad  g_j(\la):=(\la-\la_k)H_j(\la)\]
$f_j(\la)$ has a unique zero $z_j$, and we have on the circle
$C(z_j,\eps_j)$ with centre $z_j$ and radius $\eps_j$
\[|f_j(\la)|=\eps_j ,\quad \forall \la\in C(z_j,\eps_j)\]
\[|g_j(\la)|\leq (|m_k|\alpha_j^{1+2k}+\eps_j)|H_j(\la)|\leq
 C\alpha_j^\demi \eps_j, \quad \forall \la\in C(z_j,\eps_j)\]
with $C>0$ which does not depend on $j$.
We then deduce that there exists an integer $J_1\geq J_0$
such that for all $j\geq J_1$
\[|g_j(\la)|<|f_j(\la)|, \quad \forall \la\in C(z_j,\eps_j)\]
what ensures, by Rouch\'e's theorem, that $f_j+g_j$ has exactly one
zero in the disc whose boundary is $C(z_j,\eps_j)$.
But $(f_j+g_j)(\la_k)=m_k\alpha_j^{1+2k}\not=0$ thus
$S_j(\la)^{-1}$ has a unique pole in the disc whose boundary
is $C(z_j,\eps_j)$ and $S(\la)^{-1}$ has a sequence of poles converging
to $\la_k$. By the inversion formula  $S(n-\la)=S(\la)^{-1}$,
we deduce that there exists a sequence of  resonances
approaching $\frac{n-1}{2}-k$.
\qed

\end{document}